%%%%%%%%%%%%%%%%%%%%%%%%%%%%%%%%%%%%%%%%%%%%%%%%%%%%%%%%%%
%%%%%%%%%%%%%%%%%%%%%%%%%%%%%%%%%%%%%%%%%%%%%%%%%%%%%%%%%%
%%
%%     This is the AMS-LaTeX file:
%%
%%     Colli Gilardi 1
%%     Paper dedicated to Juergen 

%%
%%%%%%%%%%%%%%%%%%%%%%%%%%%%%%%%%%%%%%%%%%%%%%%%%%%%%%%%%%

\def\LaTeX{%
  \let\Begin\begin
  \let\End\end
  \let\salta\relax
  \let\finqui\relax
  \let\futuro\relax}

\def\UK{\def\our{our}\let\sz s}
\def\USA{\def\our{or}\let\sz z}

\UK
%\USA

%%%%%%%%%%%%%%%%%%%%%%%%%%%%%%%%%

% scegliere fra \TeX e \LaTeX  e fra  \UK oppure \USA

%\TeX
\LaTeX

%\UK
\USA

%%%%%%%%%%%%%%%%%%%%%%%%%%%%%%%%%
%% page layout
%%%%%%%%%%%%%%%%%%%%%%%%%%%%%%%%%

\salta

\documentclass[twoside,12pt]{article}
\setlength{\textheight}{24cm}
\setlength{\textwidth}{16cm}
\setlength{\oddsidemargin}{2mm}
\setlength{\evensidemargin}{2mm}
\setlength{\topmargin}{-15mm}
\parskip2mm

%%%%%%%%%%%%%%%%%%%%%%%%%%%%%%%%%
%% packages
%%%%%%%%%%%%%%%%%%%%%%%%%%%%%%%%%

%\usepackage{color}
\usepackage[usenames,dvipsnames]{color}
\usepackage{amsmath}
\usepackage{amsthm}
\usepackage{amssymb}
\usepackage[mathcal]{euscript}

%\usepackage[notref,notcite]{showkeys}
%\usepackage{showkeys}
%
%		COLORS FOR CORRECTIONS
%
% do the same, please (i.e., don't use the standard {\color{red} text} or similar): 
% just choose the color you prefer in \def\yourname

% example of use:  \juerg{I want this to become blue}

\definecolor{viola}{rgb}{0.3,0,0.7}
\definecolor{ciclamino}{rgb}{0.5,0,0.5}

\def\gianni #1{#1}
\def\pier #1{#1}
\def\piernew #1{#1}

%\def\revis #1{{\color{red}#1}}

%\def\pier #1{#1}
%\def\gianni #1{#1}

%%%%%%%%%%%%%%%%%%%%%%%%%%%%%%%%%
%% you may adjust the baseline
%%%%%%%%%%%%%%%%%%%%%%%%%%%%%%%%%

%\renewcommand{\baselinestretch}{0.95}

%%%%%%%%%%%%%%%%%%%%%%%%%%%%%%%%%
%% bibliographystyle
%%%%%%%%%%%%%%%%%%%%%%%%%%%%%%%%%

\bibliographystyle{plain}

%%%%%%%%%%%%%%%%%%%%%%%%%%%%%%%%%
%% environments
%%%%%%%%%%%%%%%%%%%%%%%%%%%%%%%%%

%

\finqui

\def\Beq{\Begin{equation}}
\def\Eeq{\End{equation}}
\def\Bsist{\Begin{eqnarray}}
\def\Esist{\End{eqnarray}}

\def\Bthm{\Begin{theorem}}
\def\Ethm{\End{theorem}}
\def\Blem{\Begin{lemma}}
\def\Elem{\End{lemma}}
\def\Bprop{\Begin{proposition}}
\def\Eprop{\End{proposition}}

\def\Brem{\Begin{remark}\rm}
\def\Erem{\End{remark}}

\def\Bdim{\Begin{proof}}
\def\Edim{\End{proof}}
\def\Bcenter{\Begin{center}}
\def\Ecenter{\End{center}}
\let\non\nonumber

%%%%%%%%%%%%%%%%%%%%%%%%%%%%%%%%%
%% macros
%%%%%%%%%%%%%%%%%%%%%%%%%%%%%%%%%

% macro salvate

% sottosezioni non numerate

\def\step #1 \par{\medskip\noindent{\bf #1.}\quad}

% abbreviazioni di parole

\def\Lip{Lip\-schitz}
\def\Holder{H\"older}

\def\aand{\quad\hbox{and}\quad}

\def\lhs{left-hand side}
\def\rhs{right-hand side}
\def\sfw{straightforward}
\def\omegalimit{$\omega$-limit}

% versioni inglesi (UK) o americane (USA)

\def\bhv{behavi\our}

% bold, cal e mathop

\def\multibold #1{\def\arg{#1}%
  \ifx\arg\pto \let\next\relax
  \else
  \def\next{\expandafter
    \def\csname #1#1#1\endcsname{{\bf #1}}%
    \multibold}%
  \fi \next}

\def\pto{.}

\def\multical #1{\def\arg{#1}%
  \ifx\arg\pto \let\next\relax
  \else
  \def\next{\expandafter
    \def\csname cal#1\endcsname{{\cal #1}}%
    \multical}%
  \fi \next}

% operatori

\def\multimathop #1 {\def\arg{#1}%
  \ifx\arg\pto \let\next\relax
  \else
  \def\next{\expandafter
    \def\csname #1\endcsname{\mathop{\rm #1}\nolimits}%
    \multimathop}%
  \fi \next}

\multibold
qwertyuiopasdfghjklzxcvbnmQWERTYUIOPASDFGHJKLZXCVBNM.

\multical
QWERTYUIOPASDFGHJKLZXCVBNM.

\multimathop
diag dist div dom mean meas sign supp .

% accorpamenti di formule citate:
% uso  \accorpa {prima}{seconda}
%      \Accorpa\cs prima seconda (con il comodo blank anche dopo)
% NB: \Accorpa definisce \cs come l'accorpamento delle due citazioni
% e scrive sul file.log

\def\accorpa #1#2{\eqref{#1}--\eqref{#2}}
\def\Accorpa #1#2 #3 {\gdef #1{\eqref{#2}--\eqref{#3}}%
  \wlog{}\wlog{\string #1 -> #2 - #3}\wlog{}}

% macro comode

\def\separa{\noalign{\allowbreak}}

\def\somma #1#2#3{\sum_{#1=#2}^{#3}}

\def\graffe #1{\mathopen\{#1\mathclose\}}

\def\<#1>{\mathopen\langle #1\mathclose\rangle}
\def\norma #1{\mathopen \| #1\mathclose \|}

\def\[#1]{\mathopen\langle\!\langle #1\mathclose\rangle\!\rangle}

\def\iot {\int_0^t}
\def\ioT {\int_0^T}
\def\intQt{\int_{Q_t}}
\def\intQ{\int_Q}
\def\iO{\int_\Omega}

\def\dt{\partial_t}
\def\dn{\partial_\nu}

\def\cpto{\,\cdot\,}

\def\checkmmode #1{\relax\ifmmode\hbox{#1}\else{#1}\fi}
\def\aeO{\checkmmode{a.e.\ in~$\Omega$}}
\def\aeQ{\checkmmode{a.e.\ in~$Q$}}

\def\aet{\checkmmode{a.e.\ in~$(0,T)$}}

\def\aaQ{\checkmmode{for a.a.~$(x,t)\in Q$}}

\def\aat{\checkmmode{for a.a.~$t\in(0,T)$}}

% insiemi numerici

\def\erre{{\mathbb{R}}}
\def\erren{\erre^n}

\def\enne{{\mathbb{N}}}

% spazi di funzioni a valori vettoriali su [0,T], [0,t], [0,s], [0,+\infty), [\delta,T]

% Come ricordare: in generale i simboli L H W  C da soli per gli spazi su (0,T)
% gli stessi raddoppiati per (0,+\infty)
% aggiunta di t o s al simbolo per (0,t) e (0,s)
% aggiunta di d al simbolo semplice o doppio per intervalli (\delta,T) e (\delta,+\infty)
% il simbolo C e i suoi derivati mettono le quadre anziche' le tonde

% Esempi   \L2V   \L\infty\Vp   \W{1,1}H   \C0H   \LL2V   \CC0\Vp   \Ld2V  \CCdH

\def\genspazio #1#2#3#4#5{#1^{#2}(#5,#4;#3)}
\def\spazio #1#2#3{\genspazio {#1}{#2}{#3}T0}

\def\L {\spazio L}
\def\H {\spazio H}

\def\C #1#2{C^{#1}([0,T];#2)}
\def\spazioinf #1#2#3{\genspazio {#1}{#2}{#3}{+\infty}0}
\def\LL {\spazioinf L}

% spazi di funzioni su \Omega, \Gamma, Q e \Sigma

\def\Lx #1{L^{#1}(\Omega)}
\def\Hx #1{H^{#1}(\Omega)}

\def\LQ #1{L^{#1}(Q)}

\def\Luno{\Lx 1}
\def\Ldue{\Lx 2}

\def\Hdue{\Hx 2}

% spazi di funzioni su Q e S

\def\LQ #1{L^{#1}(Q)}

% lettere greche

\let\theta\vartheta
\let\eps\varepsilon
\let\phi\varphi
\let\hat\widehat

\let\TeXchi\chi                         % new \chi, exactly on the baseline
\newbox\chibox
\setbox0 \hbox{\mathsurround0pt $\TeXchi$}
\setbox\chibox \hbox{\raise\dp0 \box 0 }
\def\chi{\copy\chibox}

% quadratino di fine dimostrazione

\def\QED{\hfill $\square$}

% abbreviazioni specifiche del lavoro

\def\VA #1{V_A^{#1}}
\def\VB #1{V_B^{#1}}
\def\VAn{V_{A,n}}
\def\VBn{V_{B,n}}
\def\VAm{V_{A,m}}
\def\VBm{V_{B,m}}

\def\Beta{\hat\beta}
\def\betaz{\beta^\circ}
\def\Pi{\hat\pi}
\def\Ell{\hat\ell}

\def\betaeps{\beta_\eps}
\def\Betaeps{\hat\beta_\eps}
\def\thetaeps{\theta_\eps}
\def\phieps{\phi_\eps}
\def\Feps{F_\eps}

\def\thetan{\theta^n}
\def\phin{\phi^n}

\def\tn{t_n}
\def\thetao{\theta_\omega}
\def\phio{\phi_\omega}
\def\thetas{\theta_s}
\def\phis{\phi_s}
\def\thetai{\theta^\infty}
\def\phii{\phi^\infty}

\def\thetaz{\theta_0}
\def\phiz{\phi_0}

%%%%%%%%%%%%%%%%%%%%%%%%%%%%%%
\Begin{document}
%%%%%%%%%%%%%%%%%%%%%%%%%%%%%%%%%

%%%%%%%%%%%%%%%%%%%%%%%%%%%%%%%%%
%% front page
%%%%%%%%%%%%%%%%%%%%%%%%%%%%%%%%%

\title{Well-posedness, regularity and \pier{asymptotic analyses}\\ 
  for \pier{a fractional} phase field system}

\author{}
\date{}
\maketitle
\Bcenter
\vskip-1cm
{\large\sc Pierluigi Colli$^{(*)}$}\\
{\normalsize e-mail: {\tt pierluigi.colli@unipv.it}}\\[.25cm]
{\large\sc Gianni Gilardi$^{(*)}$}\\
{\normalsize e-mail: {\tt gianni.gilardi@unipv.it}}\\[.45cm]
$^{(*)}$
{\small Dipartimento di Matematica ``F. Casorati'', Universit\`a di Pavia}\\
{\small and Research Associate at the IMATI -- C.N.R. Pavia}\\
{\small via Ferrata 5, 27100 Pavia, Italy}\\
[1cm]
{\it Dedicated to our friend Prof. Dr. J\"urgen Sprekels\\[.1cm]
on the occasion of his 70th birthday\\[.1cm]
with best wishes}
\Ecenter

\Begin{abstract}\noindent
\pier{This paper is concerned with a non-conserved phase field system of Caginalp type 
in which the main operators are fractional versions of
two fixed linear operators $A$ and $B$. The operators $A$ and $B$ are supposed to be 
densely defined, unbounded, self-adjoint, monotone in the Hilbert space $L^2(\Omega)$,
for some bounded and smooth domain $\Omega$, and have compact resolvents. Our 
definition of the fractional powers of operators uses the approach via spectral theory.
A nonlinearity of double-well type occurs in the phase equation and either a regular or logarithmic potential, as well as a non-differentiable potential involving an
indicator function, is admitted in our approach. We show general well-posedness and regularity results, extending the corresponding results that are known for 
the non-fractional elliptic operators with zero Neumann conditions or 
other boundary conditions like Dirichlet or Robin ones. Then, we
investigate the longtime \bhv\ of the system, by fully characterizing every element of
the $\omega$-limit as a stationary solution. In the final part of the paper we 
\gianni{study} the asymptotic behavior of the system a
as the \gianni{parameter $\sigma$ appearing in the operator $B^{2\sigma}$ that plays} in the phase equation 
decreasingly tends to zero. We can prove convergence to a phase relaxation problem at the limit, in which an additional term containing the projection of the phase variable on the kernel of $B$ appears.} 
\End{abstract}
\salta
\pagestyle{myheadings}
\newcommand\testodispari{\sc Colli \ --- \ Gilardi}
\newcommand\testopari{\sc Fractional phase field systems}
\markboth{\testodispari}{\testopari}
\finqui
%
%%%%%%%%%%%%%%%%%%%%%%%%%%%%%%%%%
%% very beginning
%%%%%%%%%%%%%%%%%%%%%%%%%%%%%%%%%

\section{Introduction}
\label{Intro}
\setcounter{equation}{0}

Let $\Omega\subset\erre^3$ denote a bounded, connected and smooth set.
We deal with the Cauchy problem for the evolutionary system
\begin{align}
  & \dt\theta + \ell(\phi)\dt\phi + A^{2r} \theta = f
  \label{Iprima}
  \\
  & \dt\phi + B^{2\sigma}\phi + F'(\phi) = \theta \, \ell(\phi)
  \label{Iseconda}
\end{align}
where $A^{2r}$ and $B^{2\sigma}$, with $r>0$ and $\sigma>0$,
denote fractional powers of the self-adjoint,
monotone and unbounded linear operators $A$ and~$B$, respectively, which are densely defined 
in $H:=\Ldue$ and are supposed to have compact resolvents. 

The above system is a generalization of the well-known phase field 
system, which models a  
phase transition process taking place in the container~$\Omega$. 
In this case,
one typically has  $A^{2r}=B^{2\sigma}=-\Delta$ with\pier{, e.g., zero Neumann boundary conditions if no flux through the boundary is assumed for both variables. About the meaning 
of variables in} \accorpa{Iprima}{Iseconda}, let us notice that
the first unknown function $\theta$ \pier{represents} the \emph{relative temperature} 
near some critical value~$\theta_c$, while $\phi$ \pier{usually denotes}
the \emph{order parameter}.
Moreover, the real function $\ell$ represents the \emph{latent heat density}
and $f$ is a source term.
Finally, $F'$~denotes the derivative of \pier{a potential~$F$, which may have a double-well shape.}

\pier{Thus, the coupled equations \eqref{Iprima}--\eqref{Iseconda} yield
a system of phase field type. From the seminal work of Caginalp and coworkers (see, e.g., \cite{caginalp, cagnish}) it became clear that phase field systems are particularly suited to 
represent the dynamics of moving interfaces arising in thermally induced 
phase transitions. Typical} and physically significant examples for the potential $F$ 
are the so-called {\em classical regular potential}, the {\em logarithmic double-well potential\/},
and the {\em double obstacle potential\/}, which are given, in this order,~by
\begin{align}
  & F_{reg}(s) := \frac 14 \, (s^2-1)^2 \,,
  \quad s \in \erre, 
  \label{regpot}
  \\
  & F_{log}(s) := \bigl( (1+s)\ln (1+s)+(1-s)\ln (1-s) \bigr) - c_1 s^2 \,,
  \quad s \in (-1,1),
  \label{logpot}
  \\[1mm]
  & F_{2obs}(s) := - c_2 s^2 
  \quad \hbox{if $|s|\leq1$}
  \aand
  F_{2obs}(s) := +\infty
  \quad \hbox{if $|s|>1$}.
  \label{obspot}
\end{align}
Here, the constants $c_i$ in \eqref{logpot} and \eqref{obspot} satisfy
$c_1>1$ and $c_2>0$, so that $f_{log}$ and $f_{2obs}$ are nonconvex.
In cases like \eqref{obspot}, one has to split $F$ into a nondifferentiable convex part~\smash{$\Beta$}
(the~indicator function of $[-1,1]$, in the present example) and a smooth 
concave perturbation~$\Pi$ ($\Pi(\gianni s) = - c_2 s^2 $, $s\in \erre$, in \eqref{obspot}).
Accordingly, one has to replace the derivative of the convex part
by the subdifferential and interpret \eqref{Iseconda} as a differential inclusion
or, equivalently, as a variational inequality involving \smash{$\Beta$} rather than its subdifferential. Actually, we will mostly do the latter in this paper. 

In fact, in the present paper 
we discuss the solvability of the initial value problem for the system \eqref{Iprima}--\eqref{Iseconda} in the framework when both $r$ and $\sigma$ are positive, by first proving a well-posedness result. This is worked out in a suitable framework, first in the case of a constant $\ell$ for general operators $A^r$ and $B^\sigma$, then for a bounded and Lipschitz continuous nonlinearity $\ell$, under some restriction on the domains of $A^r$ and $B^\sigma$; indeed, $D(A^r)$ and $ D(B^\sigma)$ have to be contained into appropriate Lebesgue spaces. 
Then, in a second part of the discussion we show some regularity results and we also 
investigate the longtime \bhv\ of the system, by fully characterizing the $\omega$-limit. In the final part of the paper we focus our attention on the analysis, which turns out to be rather delicate and not trivial at all, of the asymptotic behavior of the system \eqref{Iprima}--\eqref{Iseconda}
as the \gianni{coefficient $\sigma$} playing in \eqref{Iseconda} decreases to~$0$. We can prove convergence to a phase relaxation problem at the limit, in the special case of a constant  $\ell$ and also for a concave quadratic function~$\Pi$. However, the full set of our results is precisely described in the next section, in great detail. 

Thus, here we are dealing with fractional operators, which are nowadays a challenging subject for mathematicians: in particular, different variants of fractional operators 
may be considered and tackled. Let us mention some related contribution, starting from  
the paper \cite{Kwa}, which 
deals with several definitions of the fractional Laplacian, which is a core 
example of a class of nonlocal pseudodifferential operators appearing in various areas 
of theoretical and applied mathematics.  We quote some contributions 
by Servadei and Valdinoci:
in \cite{SV1}, a comparison is made between the spectrum of two different fractional 
Laplacian operators, of which the second one fits in our framework; in \cite{SV2} the regularity of the weak solution to the fractional Laplace 
equation is discussed; \cite{SV0} is concerned with the existence of nontrivial solutions for nonlocal semilinear Dirichlet problem; in \cite{SV3} the authors show a fractional counterpart to the well-known Brezis--Nirenberg result on the existence of nontrivial solutions to elliptic equations with critical nonlinearities. In \cite{AD} a construction of harmonic functions on bounded domains is given for the spectral fractional Laplacian operator. In the paper \cite{CT}, the authors investigate a nonlinear pseudodifferential boundary value problem 
in a bounded domain with homogeneous Dirichlet boundary conditions. 
Regularity results and sharp estimates are discussed in \cite{CS} for fractional elliptic equations. Fractional Dirichlet and Neumann type 
boundary problems associated with the fractional Laplacian are investigated in  
\cite{Gru1}, by {demonstrating} regularity properties with a spectral approach; this 
analysis is extended to the fractional heat equation in \cite{Gru2}. The contribution 
\cite{MN} deals with obstacle problems for the spectral fractional Laplacian   
The authors of \cite{RS1, RS2} prove regularity up to the boundary for a 
boundary value problem using the Caputo variant of an integral operator 
with the Riesz kernel. Some nonlocal problems involving the fractional $p -$Laplacian and nonlinearities at critical growth are examined in \cite{BSY}.
Fractional porous medium type equations are discussed in \cite{BFV,BSV,BV}:
\cite{BSV} deals with existence, uniqueness and asymptotic behavior of the 
solutions to an integro-differential equation related to porous medium equations {in}  
bounded domains; uniform 
estimates for positive solutions of a porous medium equation are derived in \cite{BV}, 
where the spectral fractional Laplacian with zero Dirichlet boundary data is 
considered; a quantitative study of nonnegative solutions of the same equation 
is provided in \cite{BFV}, where decay and positivity, Harnack inequalities, 
interior and boundary regularity, and asymptotic behavior are investigated.

We point out that there are already some contributions addressing 
nonlocal variants of Allen-Cahn, Cahn--Hilliard and phase field 
systems. In \cite{AkSS1}, Akagi, Schimperna and Segatti introduce 
a fractional variant of the Cahn--Hilliard equation settled 
in a bounded domain and complemented with homogeneous Dirichlet boundary conditions of 
solid type; they prove existence and uniqueness of weak solutions and 
investigate some significant singular limits as the order of either of the fractional 
Laplacians appearing in the system approaches zero. In this respect, their results can be compared with our results of Section~\ref{RELAXATION}: it is worth mentioning that the authors of  \cite{AkSS1}
use fractional operators not defined via the spectal properties and actually different 
from ours. In the recent paper \cite{AkSS2}, for fixed 
orders of the operators, the same authors show the convergence as time goes to infinity of each solution to a (single) equilibrium. The contribution \cite{AM} deals with a fractional Cahn--Hilliard 
equation by considering a gradient flow in the negative order Sobolev space $H^{-\alpha}$, 
$\alpha\in [0,1]$, 
where the case $\alpha=1$ corresponds to the classical Cahn--Hilliard equation and 
the choice $\alpha=0$ recovers the Allen--Cahn equation; existence and stability estimates 
are proved. We also mention the articles \cite{g,gz} that are concerned with nonlocal
phase field models for phase separations, using a free energy which arises naturally in 
the analysis of the large scale limit of systems of interacting particles. Another interesting 
analysis 
of a nonstandard and nonlocal Cahn--Hilliard system can be found in \cite{CGS6}. A non-local 
version of the Cahn--Hilliard equation is treated in \cite{GalEJAM}; the papers 
\cite{GalDCDS,GalAIHP} investigate a doubly nonlocal Cahn--Hilliard equation with special 
kernels in the operators, by focusing on the interaction between the two levels 
of nonlocality in the corresponding terms. The paper \cite{HTY} studies numerical solutions to the Allen-Cahn equation with a fractional
Laplacian: the authors use the second-order Crank-Nicolson scheme to discretize the equation in time and the second-order central difference scheme for discretization in space.
A space-time fractional Allen-Cahn phase-field model that describes the transport of the fluid mixture of two immiscible fluid phases is discussed in \cite{LWY}; the space and time fractional order parameters control the sharpness and the decay behavior of the interface.
We  also quote the contribution \cite{CCFT}, where melting and solidification for metallurgical processes concerned  with phase transitions of pure metals are studied; during the solid phase 
the metals show an evident ductility and these particular phenomena can be well described by a phase field fractional model, whose evolution has to satisfy a Ginzburg-Landau equation.

In our approach,  which follows closely the setting 
recently used in \cite{CGS18}, we work with fractional 
operators defined via spectral theory. By this, we can easily consider
powers of a second-order elliptic operator with either Dirichlet or Neumann or 
Robin boundary conditions, as well as other operators, e.g., 
fourth-order ones or systems involving the Stokes operator.
The contents of the paper can be summarized here. In Section~\ref{STATEMENT}, a precise 
statement of the problem along with a full set of assumptions is given and 
most of the results proved in the paper are stated. 
Section~\ref{UNIQUENESS} deals with the continuous 
dependence of the solution on the data, while Section~\ref{APPROX} introduces 
an approximating problem based on the Moreau--Yosida 
regularizations of the convex function and on a 
Faedo-Galerkin scheme, which is sharply discussed about 
existence of the approximating solution and the proof of a priori estimates.
In Section~\ref{EXISTENCE} the existence proof is terminated,  
by taking the limits with respect to the Yosida approximation parameters, 
and the further regularity results are derived.  
Section~\ref{LONGTIME} brings the analysis of the long-time \bhv\ and the 
characterization of the $\omega$-limit as set of stationary solutions to the system 
\eqref{Iprima}--\eqref{Iseconda}. Finally, 
Section~\ref{RELAXATION} is completely devoted to the study of 
the asymptotic behavior of the system \eqref{Iprima}--\eqref{Iseconda}
as the \gianni{parameter $\sigma$} tends to~$0$: the 
convergence to a phase relaxation problem at the limit is rigorously proved. 

%%%%%%%%%%%%%%%%%%%%%%%%%%%%%%%%%%%%%%%%%%%%%%%%%%%%%%%%%%%%%%%%%%%%%%%%

\section{Statement of the problem and results}
\label{STATEMENT}
\setcounter{equation}{0}

In this section, we state precise assumptions and notations and present our results.
First of all, the set $\Omega\subset\erre^3$ is assumed to be bounded, connected and 
smooth, and $\nu$ and $\dn$ denote the outward unit normal vector field on $\Gamma:=\partial\Omega$
and the corresponding normal derivative, respectively.
In order to simplify the notation, we~set
\Beq
  H := \Ldue
  \label{defH}
\Eeq
and endow $H$ with its standard norm $\norma\cpto$ and inner product $(\cpto,\cpto)$.
As far as our assumptions are concerned, we first postulate that
\begin{align}
  & A:D(A)\subset H\to H
  \aand
  B:D(B)\subset H\to H
  \quad \hbox{are}
  \non
  \\
  & \hbox{unbounded monotone self-adjoint linear operators with compact resolvents} 
  \qquad
  \label{hpAB} 
\end{align}
and introduce some function spaces and fractional operators.
We avoid the \pier{background} of interpolation theory and give direct definitions.
The above assumption implies that there are sequences 
$\{\lambda_j\}$ and $\{\mu_j\}$ of eigenvalues
and orthonormal sequences $\{e_j\}$ and $\{\eta_j\}$ of corresponding eigenvectors,
that~is,
\Beq
  A e_j = \lambda_j e_j, \quad
  B \eta_j = \mu_j \eta_j
  \aand
  (e_i,e_j) = (\eta_i,\eta_j) = \delta_{ij}
  \quad \hbox{for $i,j=1,2,\dots$}
  \label{eigen}
\Eeq
such that
\begin{align}
  & 0 \leq \lambda_1 \leq \lambda_2 \leq \dots
  \aand
  0 \leq \mu_1 \leq \mu_2 \leq \dots
  \non
  \\
  & \quad \hbox{with} \quad
  \lim_{j\to\infty} \lambda_j
  = \lim_{j\to\infty} \mu_j
  = + \infty,
  \label{eigenvalues}
  \\[1mm]
  & \hbox{$\{e_j\}$ and $\{\eta_j\}$ are complete systems in $H$}.
  \label{complete}
\end{align}
\pier{These} assumptions allow us 
to introduce the Hilbert spaces $\VA r$ and $\VB\sigma$
and the power operators $A^r$ and~$B^\sigma$ 
(for \pier{some} arbitrary positive real exponents \pier{$r$ and $\sigma$}) 
as~follows
\begin{align}
  & \VA r := D(A^r)
  = \Bigl\{ v\in H:\ \somma j1\infty |\lambda_j^r (v,e_j)|^2 < +\infty \Bigr\},
  \label{defdomAr}
  \\
  & \VB\sigma := D(B^\sigma)
  = \Bigl\{ v\in H:\ \somma j1\infty |\mu_j^\sigma (v,\eta_j)|^2 < +\infty \Bigr\},
  \label{defdomBs}
  \\
  & A^r v = \somma j1\infty \lambda_j^r (v,e_j) e_j
  \aand
  B^\sigma v = \somma j1\infty \mu_j^\sigma (v,\eta_j) \eta_j
  \non
  \\
  & \quad \hbox{for $v\in\VA r$ and $v\in\VB\sigma$, respectively.}
  \label{defArBs}
\end{align}
\pier{Note that the series in \eqref{defArBs} are} convergent in the strong topology of~$H$
due to the properties of the coefficients.
We endow $\VA r$ and $\VB\sigma$ with the graph norms and inner products
\begin{align}
  &\norma v_{A,r}^2 := (v,v)_{A,r}
  \aand
  (v,w)_{A,r} := (v,w) + (A^r v , A^r w)
  \label{defnormaAr}
  \\
  &\norma v_{B,\sigma}^2 := (v,v)_{B,\sigma}
  \aand
  (v,w)_{B,\sigma} := (v,w) + (B^\sigma v , B^\sigma w)
  \label{defnormaBs}
\end{align}
for $v,w\in\VA r$ and $v,w\in\VB\sigma$, respectively.
If $r_i$ and $\sigma_i$ are arbitrary positive exponents,
it is clear that
\begin{align}
  & (A^{r_1+r_2} v,w)
  = (A^{r_1} v, A^{r_2} w)
  \quad \hbox{for every $v\in\VA{r_1+r_2}$ and $w\in\VA{r_2}$},
  \label{propA}
  \\[1mm]
  & (B^{\sigma_1+\sigma_2} v,w)
  = (B^{\sigma_1} v, B^{\sigma_2} w)
  \quad \hbox{for every  $v\in\VB{\sigma_1+\sigma_2}$ and $w\in\VB{\sigma_2}$}.
  \label{propB}
\end{align}
Moreover, since $A^r$ and $B^\sigma$ are symmetric, \pier{for $r,\sigma>0$ we also have} 
\Beq
  (\dt v,A^{2r}v) = \frac 12 \, \frac d{dt} \, \norma{A^rv}
  \aand
  \pier{(\dt w,B^{2\sigma}w) = \frac 12 \, \frac d{dt} \, \norma{B^\sigma w}}
  \label{propAB}
\Eeq
for every $v\in\H1H\cap{\L2{\VA{2r}}}$ and $\pier{w} \in\H1H\cap\L2{\VB{2\sigma}}$, respectively.
We also remark that for every $r,\sigma>0$
\Beq
  \hbox{the embeddings $\VA r\subset H$ and $\VB\sigma\subset H$ are compact},
  \label{compact}
\Eeq
as one immediately sees by using \eqref{eigenvalues}.

\pier{\Brem
\label{Remhpsimple}
Let us mention some simple situation for possible operators $A$ and $B$. 
In view of \eqref{eigenvalues} and \eqref{complete}, a standard 
elliptic \gianni{operator} with Dirichlet boundary conditions provides 
an example with a strictly positive first eigenvalue. Another 
operator that can be considered is the Laplace 
operator $-\Delta$ with Neumann boundary conditions,
which corresponds to the choice $D(-\Delta)=\{v\in\Hdue:\ \dn v=0\}$;
in this case the first eigenvalue is $0$ and it is simple, with 
corresponding eigenfunctions that are the constant functions.
\Erem}

Coming back to our system, we fix $r$ and $\sigma$ once and for all.
Thus, we postulate that
\Beq
  r, \sigma \in (0,+\infty).
  \label{hpconst}
\Eeq
For the nonlinearities, we require the properties listed below
and notice that they are fulfilled by all of the \pier{significant} potentials 
\accorpa{regpot}{obspot}.
We assume that
\Bsist
  & \Beta : \erre \to [0,+\infty]
  \quad \hbox{is convex, proper and l.s.c.\ with} \quad
  \Beta(0) = 0,
  \label{hpBeta}
  \\
  \separa
  & \Pi : \erre \to \erre
  \quad \hbox{is of class $C^1$ with a \Lip\ continuous first derivative}.
  \label{hpPi}
\Esist
We set, for convenience,
\Beq
  \beta := \partial\Beta 
  \aand
  \pi := \Pi' .
  \label{defbetapi}  
\Eeq
\Accorpa\HPstruttura hpconst defbetapi
Moreover, we term $D(\Beta)$ and $D(\beta)$ the effective domains of~$\Beta$ and~$\beta$, respectively,
and, for $r\in D(\beta)$,
we use the symbol $\betaz(r)$ for the element of $\beta(r)$ having minimum modulus.
We notice that $\beta$ is a maximal monotone graph in $\erre\times\erre$.
We also remark that \eqref{hpPi} implies that 
$\Pi$ grows at most quadratically and that $\pi$ is linearly bounded.
Finally, we assume that
\Beq
  \ell : \erre \to \erre
  \quad \hbox{is bounded and \Lip\ continuous}.
  \label{hpell}
\Eeq
However, in order to keep the \pier{operators} $A^r$ and $B^\sigma$ as general as possible,
we often assume that $\ell$ is a constant.
Indeed, the more general case~\eqref{hpell} needs further assumptions.

At this point, we can state the problem we aim to discuss.
While the first equation coincides with~\eqref{Iprima},
we present \eqref{Iseconda} as a variational inequality
written  in a weak form on account of \eqref{propB}.
For the data, we make the following assumptions:
\begin{align}
  & f \in \L2H 
  \label{hpf}
  \\
  & \thetaz \in \VA r , \quad
  \phiz \in \VB\sigma
  \aand
  \Beta(\phiz) \in \Luno .
  \label{hpz}
\end{align}
\Accorpa\HPdati hpf hpz
Then, we set
\Beq
  Q := \Omega \times (0,T)
  \label{defQ}
\Eeq
and look for a pair $(\theta,\phi)$ satisfying
\begin{align}
  & \theta \in \H1H \cap \L\infty{\VA r} \cap \L2{\VA{2r}},
  \qquad
  \label{regtheta}
  \\
  & \phi \in \H1H \cap \L\infty{\VB\sigma} ,
  \label{regphi}
  \\
  & \Beta(\phi) \in \LQ1
  \label{regBetaphi}
\end{align}
\Accorpa\Regsoluz regtheta regBetaphi
and solving the system
\begin{align}
  & \dt\theta + \ell(\phi) \dt\phi + A^{2r} \theta = f
  \quad \aeQ \, ,
  \label{prima}
  \\[1mm]
  & \bigl( \dt\phi(t) , \phi(t) - v \bigr)
  + \bigl( B^\sigma \phi(t) , B^\sigma (\phi(t) - v) \bigr) 
%  \non
%  \\
%  & \quad {}
  + \iO \Beta(\phi(t))
  + \bigl( \pi(\phi(t)) , \phi(t) - v \bigr)
  \non
  \\
  & \leq \bigl( \theta(t) \, \ell(\phi(t)) , \phi(t) - v \bigr)
  + \iO \Beta(v)
%  \non
%  \\
%  & 
  \quad 
  \hbox{\aat\ and every $v\in\VB\sigma$} \, ,
  \label{seconda}
  \\[1mm]
  & \theta(0) = \thetaz
  \aand
  \phi(0) = \phiz \,.
  \label{cauchy}
\end{align}
\Accorpa\Pbl prima cauchy
We notice that \pier{equation~\eqref{prima} has been written a.e.\ in $Q$ and in this case the single terms, including 
$A^{2r} \theta$, are interpreted as functions of space and time; another way of reading~\eqref{prima} could be 
a.e.\ in $(0,T)$ as the equality makes sense for all terms in the space $H$ as well. In our notation, here and in the sequel, we follow the former approach. Concerning \eqref{seconda}, we warn the reader that 
$$ \iO \Beta(v) = +\infty  \quad \hbox{whenever} \quad \Beta (v) \notin \Luno.$$
We follow a similar agreement for integrals of the type $\intQ \Beta(v)$ whenever 
$v\in L^2(Q)$ but $\Beta(v) \notin L^1(Q)$. We also remark that}
\eqref{seconda} is equivalent to the following time-integrated \pier{version:}
\begin{align}
  & \intQ \dt\phi \, (\phi-v)
  + \intQ B^\sigma \phi \, B^\sigma(\phi-v)
  + \intQ \Beta(\phi)
  + \intQ \pi(\phi) \, (\phi-v)
  \non
  \\
  & \leq \intQ \theta \, \ell(\phi) (\phi-v)
  + \intQ \Beta(v)
  \qquad \hbox{for every $v\in\L2{\VB\sigma}$}.
  \label{intseconda}
\end{align}
\pier{\Brem
\label{Signif}
According to the definition of subdifferential (cf., e.g., \cite{Brezis} or \cite{Barbu}), the precise meaning of the inequality \eqref{seconda} is that there exists some element $\xi \in \L2{(V_B^\sigma)^*}$ such that 
$$ \xi:= \theta \, \ell(\phi) - \dt \phi -  B^{2\sigma} \phi -  \pi(\phi) \in \partial \Phi  (y)   \quad \hbox{a.e. in } (0,T) ,$$
where $\partial \Phi $ is the subdifferential of the convex function 
$\Phi : V_B^\sigma \to [0, +\infty]$ defined by 
$$
\Phi (v) :=  \iO \Beta(v) \quad \hbox{if} \ \, \Beta(v) \in \Luno, \quad \Phi (v) := +\infty \quad \hbox{otherwise.}
$$
Indeed, we point out that the subdifferential $\partial \Phi $ is a maximal monotone operator 
from $V_B^\sigma$ to $(V_B^\sigma)^*$. In this sense, \eqref{seconda} turns out to be 
a slight generalization of \eqref{Iseconda}. 
\Erem}

\pier{The assumptions} \HPstruttura\ we have made till now on the structure are very general.
Nevertheless, they are sufficient to guarantee well-posedness and continuous dependence
at least if $\ell$ is linear
(the nonlinear case is discussed later~on).
In the result stated below and in the rest of the paper,
for $v\in\LQ1$, the symbol $1*v$ denotes the function belonging to $\LQ1$ that is defined~by
\Beq
  (1*v)(x,t) := \iot v(x,s) \, ds
  \quad \aaQ .
  \label{convoluz}
\Eeq

\Bthm
\label{Wellposedness}
Assume that \HPstruttura\ are satisfied and that $\ell$ is a \pier{constant}.
Moreover, let the assumptions \HPdati\ on the data be fulfilled.
Then there exists a unique pair $(\theta,\phi)$ 
satisfying \Regsoluz\ and solving problem \Pbl. 
Moreover, if $(f_i,\theta_{0i},\phi_{0i})$, $i=1,2$, are two choices of the data 
and $(\theta_i,\phi_i)$ are the corresponding solutions, then
we have
\begin{align}
  & \norma{\theta_1-\theta_2}_{\L2H}
  + \norma{1*(\theta_1-\theta_2)}_{\L\infty{\VA r}}
  + \norma{\phi_1-\phi_2}_{\L\infty H\cap\L2{\VB\sigma}}
  \non
  \\
  & \leq K \bigl(
    \norma{1*(f_1-f_2)}_{\L2H}
    + \norma{\theta_{01}-\theta_{02}}
    + \norma{\phi_{01}-\phi_{02}}
  \bigr)
  \label{contdep}
\end{align}
with a constant $K$ that depends only on \pier{$\ell$, some Lipschitz constant for $\pi$,} and~$T$.
\Ethm

At this point, one can wonder whether both $B^{2\sigma}\phi$ and $\beta(\phi)$ make sense in $\LQ2$
and \eqref{seconda} yields something that is closer to~\eqref{Iseconda}, like
\Beq
  \dt\phi + B^{2\sigma}\phi + \xi + \pi(\phi) = \ell\theta
  \quad \aeQ 
  \label{secondaqo}
\Eeq
for some function $\xi$ on $Q$ satisfying $\xi\in\beta(\phi)$ \aeQ.
This depends on the assumption
\Beq
  \betaeps(v) \in \VB\sigma
  \aand
  \bigl( B^\sigma \betaeps(v) , B^\sigma v \bigr) \geq 0
  \quad \hbox{for every $v\in\VB\sigma$ and $\eps>0$}
  \label{hpqo}
\Eeq
where $\betaeps$ denotes the Yosida regularization of $\beta$ at the level~$\eps>0$
(see, e.g., \cite[p.~28]{Brezis}).
We notice that \eqref{hpqo} does not follow from \HPstruttura\ as a consequence
and is rather restrictive.
Essentially, it is fulfilled \pier{whenever} $B^{2\sigma}$ is one of the more usual second order linear elliptic operators
with boundary conditions of a standard type, indeed. 
Therefore, in the general case of fractional powers, 
the more proper formulation of the equation \eqref{Iseconda} for $\phi$ 
is the variational inequality~\eqref{seconda} \pier{(see also Remark~\ref{Signif})}. 

\Bprop
\label{Secondaqo}
In addition to the assumptions of Theorem~\ref{Wellposedness},
suppose that \eqref{hpqo} is fulfilled and let $(\theta,\phi)$ be the solution to \Pbl.
Then, $\phi$ enjoys the regularity property
\Beq
  \phi \in \L2{\VB{2\sigma}}
  \label{piuregphi}
\Eeq
and there exists $\xi$ satisfying
\Beq
  \xi \in \L2H
  \aand
  \xi \in \beta(\phi)
  \quad \aeQ\, ,
  \label{regxi}
\Eeq
such that the differential equation \eqref{secondaqo} holds true.
\Eprop

Independently of assumption \eqref{hpqo} and of the above result,
we can show some more regularity for the solution if the datum $\phiz$ satisfies some proper conditions,
as stated below.

\Bthm
\label{Regularity}
In addition to the assumptions of Theorem~\ref{Wellposedness}, suppose that
\Beq
  \phiz \in \VB{2\sigma}
  \aand
  \betaz(\phiz) \in H .
  \label{hpregdtphi}
\Eeq
Then, the solution $(\theta,\phi)$ to problem \Pbl\ also satisfies
\Beq
  \dt\phi \in \L\infty H \cap \L2{\VB\sigma}.
  \label{regdtphi}
\Eeq
%\pier{NON FUNZIONA: If in addition $f\in\L\infty H$, then}
%\Beq
%  \pier{\dt\theta \in \L\infty H
%  \aand
%  \theta \in \L\infty{\VA{2r}}.}
%  \label{regdttheta}
%\Eeq
\Ethm
\Brem
\label{Estimates}
Of course, to each of our existence and regularity results
one can associate a bound for some norm of the solution
through a constant that depends only on the assumptions at hand and~$T$.
Such bounds are obtained from the construction of the solution, directly.
For instance, concerning Theorem~\ref{Wellposedness},
we have the following estimate 
\Beq
  \norma\theta_{\H1H\cap\L\infty{\VA r}\cap\L2{\VA{2r}}}
  + \norma\phi_{\H1H\cap \L\infty{\VB\sigma}}
  \leq K_1\, , 
  \label{stimasoluz}
\Eeq
where the constant $K_1$ depends only on the structure of the system, 
the norms of the data corresponding to \HPdati, and~$T$.
\Erem

\pier{Concerning the case of a non-constant $\ell$ that 
satisfies just~\eqref{hpell}, we can show some results which}
depend on further assumptions on the operators $A^r$ and~$B^\sigma$.
Namely, we require the following Sobolev-type \pier{embeddings}:
\begin{align}
  & \hbox{there exist $p,q\in[1,+\infty]$ with} \quad
  \frac 1p + \frac 2q = 1
  \quad \hbox{such that} 
  \non \ 
%  \\
%  & 
  \VA r \subset \Lx p
%  \\
 \ 
  \hbox{and}
  \\
  &\VB\sigma \subset \Lx q ,
%  \non 
  \ 
%  \\
%  & 
  \hbox{the embeddings being continuous and compact, respectively}.
  \label{hpsobolev}
\end{align}
We notice that this implies the compactness inequality
\Beq
  \norma v_q \leq \delta \, \norma{B^\sigma v} + C_\delta \, \norma v
  \quad \hbox{for every $v\in\VB\sigma$ and $\delta>0$},
  \label{compactineq}
\Eeq
with some constant $c_\delta$ depending on~$\delta$.
Of course, $\norma\cpto_q$ denotes the norm in~$\Lx q$.
This notation is also used in the next sections.

\Brem
\label{Sobolev}
It is worth noting that assumption \eqref{hpsobolev}
is satisfied in the case of standard second order operators $A$ and $B$
provided that $r$ and $\sigma$ are not too small.
Assume, for instance, that $A$ and $B$ 
are the Laplace operators with either Dirichelt or Neumann boundary conditions.
Then, for $r,\sigma\in(0,1)$, the spaces $\VA r$ and $\VB\sigma$
are embedded into the fractional Sobolev spaces $\Hx{2r}$ and~$\Hx{2\sigma}$, respectively.
Therefore, by also assuming $r,\sigma<3/4$ for simplicity,
we see that the three-dimensional embeddings required in \eqref{hpsobolev}
hold true provided that
\Beq
  2r - \frac 32 \geq - \frac 3p
  \aand
  2\sigma - \frac 32 > - \frac 3q \,,
  \quad \hbox{or} \quad
  \frac 1p \geq \frac 12 - \frac {2r}3
  \aand
  \frac 2q > 1 - \frac {4\sigma}3 \,.
  \non
\Eeq
Hence, the existence of $p$ and $q$ as in \eqref{hpsobolev} is ensured if in addition it is assumed that
\Beq
  \Bigl( \frac 12 - \frac {2r}3 \Bigr)
  + \Bigl( 1 - \frac {4\sigma}3 \Bigr)
  < 1 \,,
  \quad \hbox{i.e.,} \quad
  r + 2\sigma > \frac 34 \,.
  \non
\Eeq
\Erem

\Bthm
\label{Nonlinearell}
Besides \HPstruttura, assume that \eqref{hpell} and \eqref{hpsobolev} are fulfilled.
Moreover, \pier{let the assumptions \HPdati\ on the data be satisfied}.
Then, the same conclusions of Theorem~\ref{Wellposedness} hold true\pier{, but with a constant $K$
also depending on the norms of the data $(f_i, \theta_{0i},  \phi_{0i})$, $i=1,2$, in \HPdati.}
\Ethm

\Bprop
\label{Nonlinearsecondaqo}
In addition to the assumptions of Theorem~\ref{Nonlinearell},
suppose that \eqref{hpqo} is fulfilled.
Then, the same conclusions of Proposition~\ref{Secondaqo} hold true.
\Eprop

\Bthm
\label{Nonlinearregularity}
In addition to the assumptions of Theorem~\ref{Nonlinearell}, 
assume~\eqref{hpregdtphi}.
Then, the same conclusions of Theorem~\ref{Regularity} hold true.
\Ethm

The \pier{subsequent} aim of this paper is the study \pier{of} the longtime \bhv\ of the solution.
Precisely, under the structural assumptions 
postulated in one of Theorems~\ref{Wellposedness} and~\ref{Nonlinearell}
and the assumptions \eqref{hpz} on the initial data,
if $f$ is defined on the whole half-line $t\geq0$
and satisfies \eqref{hpf} for every $T>0$,
the existence of a unique global solution defined in $[0,+\infty)$
and satisfying \Regsoluz\ for every $T>0$ is guaranteed.
However, in order to treat \pier{its longtime} \bhv,
we need further assumptions that do not follow from the other hypotheses
and are commented in the next remarks.
First of all, we postulate the following coercivity condition
\begin{align}
  & \hbox{there exist \pier{some} positive constants $\alpha$ and $C$ such~that}
  \non
  \\
  & \Betaeps(r) + \Pi(r) \geq \alpha \, r^2 - C
  \quad \hbox{for every $r\in\erre$ and \pier{for} $\eps>0$ small enough,}
  \label{hpcoerc}
\end{align}
where $\Betaeps$ is the Moreau regularization of $\Beta$ at the level~$\eps$
(see, e.g., \cite[Prop.~2.11 p.~39]{Brezis}).
If $\ell$ is a constant, we do not need anything else on the structure of the system.
On the contrary, if $\ell$ just satisfies~\eqref{hpell},
we have to reinforce the condition \eqref{hpsobolev} on the operators $A^r$ and $B^\sigma$
by requiring that
\begin{align}
  & \hbox{there exist $p,q\in[1,+\infty]$ with} \quad
  \frac 1p + \frac 1q = \frac 12
  \quad \hbox{such that} 
  \non
 % \\ & 
  \ \, \VA r \subset \Lx p
  \, \hbox{ and}
  \\
&  \VB\sigma \subset \Lx q , 
  \hbox{the embeddings being continuous and compact, respectively.}
  \label{hpsobolevbis}
\end{align}

\Brem
\label{Coerc}
We notice that \eqref{hpcoerc} is satisfied by all \pier{of the examples} \accorpa{regpot}{obspot}.
More generally, if $\Pi(s)=-C_0\,s^2$ with $C_0>0$
(up to an additive constant) like in the quoted examples,
in order that \eqref{hpcoerc} holds true 
it is sufficient to assume that
$\Beta(s)+\Pi(s)\geq 2\alpha\,s^2 - C$ with the same $\alpha$ and~$C$.
Indeed, we have for every $s\in\erre$
\gianni{%
\begin{align}
  & \Betaeps(s)
  = \min_{\tau\in\erre} \Bigl( \frac 1{2\eps}\,(\tau-s)^2 + \Beta(\tau) \Bigr)
  \non
  \\
  & \geq \min_{\tau\in\erre} \Bigl( \frac 1{2\eps}\,(\tau-s)^2 + (2\alpha+C_0) \tau^2 - C \Bigr)
  = \bigl( 2\alpha + C_0 + O(\eps) \bigr) s^2 - C \piernew{,}
  \non
\end{align}
whence
\Beq
  \Betaeps(s) + \Pi(s)
  \geq \bigl( 2\alpha + O(\eps) \bigr) s^2 - C 
  \geq \alpha \, s^2 - C 
  \non
\Eeq
}%
if $\eps$ is small enough.

We also remark that \eqref{hpcoerc} can be weakened by replacing \gianni{$s^2$ with $|s|$} on the \rhs\
under proper assumptions on the operator~$B$ that ensure a Poincar\'e type inequality
(see \cite[Prop.~3.1]{CGS18} for a similar situation regarding the operator~$A$).
However, we assume \eqref{hpcoerc} in order to keep the linear operators as general as possible.
\Erem

\Brem
\label{LikeSobolev}
In the same framework of Remark~\ref{Sobolev}, \eqref{hpsobolevbis} is satisfied if $r+\sigma>3/4$.
\Erem

\Brem
\label{Sobolevbis}
We show that \eqref{hpsobolevbis} actually is a \pier{reinforcement} of~\eqref{hpsobolev}, i.e.,
that the former implies the latter.
Given a choice $(p_0,q_0)$ of $(p,q)$ satisfying \eqref{hpsobolevbis},
we construct $(p,q)$ fulfilling \eqref{hpsobolev}.
We take $q=q_0$, so that $\VB\sigma$ is compactly embedded in~$\Lx q$;
then, we observe that $q_0\geq2$ and define $p\in[1,+\infty]$ by means of the equality $1/p=1-(2/q_0)$,
so that $(1/p)+(2/q)=1$.
Moreover, we have that
\Beq
  \frac 1p - \frac 1{p_0}
  = \Bigl( 1 - \frac 2{q_0} \Bigr) - \Bigl( \frac 12 - \frac 1{q_0} \Bigr)
  = \frac 12 - \frac 1{q_0}
  \geq 0
  \non
\Eeq
whence $p\leq p_0$.
Hence the continuous embedding $\VA r\subset\Lx{p_0}$ we are assuming
implies the continuous embedding $\VA r\subset\Lx p$.
\Erem

Concerning the source term~$f$, we require that it tends to zero in a weak sense as time tends to infinity.
Namely, we assume~that
\Beq
  f \in \LL1H \cap \LL2H .
  \label{hpfbis}
\Eeq
Under these assumptions, we study the \omegalimit\ of $(\theta,\phi)$ in the weak topology of~$H\times H$.
This is defined as follows:
\begin{align}
  & \omega := 
  \bigl\{ (\thetao,\phio) \in H\times H : \quad
    \hbox{there esists $\{\tn\}\mathrel{\scriptstyle\nearrow}+\infty$ such that}
  \non
  \\
  & \phantom{\omega := {}} \enskip 
    (\theta(\tn),\phi(\tn)) \to (\thetao,\phio) \quad
    \hbox{weakly in $H\times H$}
  \bigr\} \,.
  \label{defomegalimit}
\end{align}
We notice that the above definition is meaningful since both $\theta$ and $\phi$ are $H$-valued continuous functions.
However, the \omegalimit\ might be empty.
Our results states that this is not the case
and that every element of $\omega$ is a pair $(\thetas,\phis)$ satisfying
\Beq
  \thetas \in \VA r
  \aand
  \phis \in \VB\sigma
  \label{regsoluzs}
\Eeq
and solving the problem
\begin{align}
  & A^r \thetas = 0
  \quad \pier{\aeO}
  \label{primas}
  \\
  & \bigl( B^\sigma\phis , B^\sigma(\phis-v) \bigr)
  + \iO \Beta(\phis)
  + \bigl( \pi(\phis) , \phis-v \bigr)
  \non
  \\
  & \leq \bigr( \thetas \, \ell(\phis) , \phis-v \bigr)
  + \iO \Beta(v)
  \qquad \hbox{for every $v\in\VB\sigma$}.
  \label{secondas}
\end{align}
We note that this system simply means that $(\thetas,\phis)$ 
is a stationary solution to problem \accorpa{prima}{seconda},
since, given any $r_1>0$, in particular $r_1=2r$, 
\eqref{primas}~is equivalent to $A^{r_1}\thetas=0$.
Indeed, such equations respectively mean the conditions
$\lambda_j^r(\thetas,e_j)=0$ for every~$j$ and $\lambda_j^{r_1}(\thetas,e_j)=0$ for every~$j$,
and the latter conditions are equivalent to each other.
However, we keep \eqref{primas} in that form for convenience.
The result we state \pier{covers} both cases regarding~$\ell$.

\Bthm
\label{Longtime}
Assume \HPstruttura, \eqref{hpcoerc} 
and either that $\ell$ is a constant or that \eqref{hpsobolevbis} and \eqref{hpell} are fulfilled.
Moreover, assume \eqref{hpz} on the initial data and \eqref{hpfbis} on the source term,
and let $(\theta,\phi)$ satisfy \pier{\eqref{regtheta}--\eqref{cauchy}} for every $T>0$.
Then, the \omegalimit\ \eqref{defomegalimit} is nonempty
and every element of it is a pair $(\thetas,\phis)$ satisfying \accorpa{regsoluzs}{secondas},
that is, it is a stationary solution.
\Ethm

\pier{The last set of results is concerned with the asymptotic \bhv\ of our system \accorpa{prima}{cauchy} as the coefficient $\sigma$ of the operator $B^\sigma$ playing in \eqref{seconda} decreases to $0$, with the aim of deducing a phase relaxation problem in the limit. 
These results are obtained in a special situation concerning the data, that is, with $\ell$ constant and also for a particular choice of the function $\pi$ (linear case). However, since 
we recognize that the present section is already rather long and in this setting we need to change a bit the notation for solutions, we prefer to postpone not only the proofs but also the statements for this part of the theory at the last section.}

The remainder of the paper is organized as follows.
The uniqueness and continuous dependence result is proved in Section~\ref{UNIQUENESS},
while the existence of a solution and its regularity are shown
in Section~\ref{EXISTENCE} and are prepared by the study 
of the approximating problem introduced in Section~\ref{APPROX}. \pier{Section~\ref{LONGTIME}
is devoted to the longtime \bhv\ of the solution. Finally, Section~\ref{RELAXATION} is concerned with the study of the limiting problem as the exponent $\sigma$ of the operator $B^\sigma$ tends to $0$.}

Throughout the paper, we widely use the notation
\Beq
  Q_t := \Omega \times (0,t)
  \quad \hbox{for $t\in(0,T]$}, \pier{\quad \hbox{with }\, \pier{Q:= Q_T},}
  \label{defQt}
\Eeq
as well as the \Holder\ inequality and the \pier{elementary} Young inequality
\Beq
  ab \leq \delta a^2 + \frac 1{4\delta} \, b^2
  \quad \hbox{for every $a,b\geq0$ and $\delta>0$}
  \label{young}
\Eeq
and we follow the general rule we explain at once concerning the constants.
The small-case italic $c$ without subscripts stands
for possibly different constants that may only depend on 
the operators $A^r$ and~$B^\sigma$,  
the shape of the nonlinearities $\beta$, $\pi$ and~$\ell$,
the properties of the data involved in the statements at hand,
and the final time~$T$,
unless some warning is given in the opposite direction.
Thus, the values of such constants do not depend on further parameters
(like the regularization parameter $\eps$ we introduce in Section~\ref{APPROX}),
and it is clear that they might change from line to line 
and even in the same formula or chain of inequalities. 
If $\pier{\delta}$ is any parameter (e.g.,~$\eps$), the symbol $c_{\pier{\delta}}$ stands for 
(possibly different) constants that depend on~$\pier{\delta}$, in addition.
In contrast, we use other symbols (e.g.,~capital letters)
for precise values of constants we want to refer~to.

%%%%%%%%%%%%%%%%%%%%%%%%%%%%%%%%%%%%%%%%%%%%%%%%%%%%%%%%%%%%%%%%%%%%%%%%

\section{Uniqueness and continuous dependence}
\label{UNIQUENESS}
\setcounter{equation}{0}

In this section, we prove the uniqueness part of Theorem~\ref{Wellposedness} 
and the continuous dependence estimate.
Moreover, we sketch how to modify our argument
for the case considered in Theorem~\ref{Nonlinearell}.
By noticing that uniqueness follows from~\eqref{contdep}
provided that this is shown for every pair of solutions,
we prove just the latter.
We fix a pair of data as in the statement and any pair of corresponding solutions
and set for convenience
\Bsist
  & f := f_1 - f_2 \,, \quad
  \thetaz := \theta_{01} - \theta_{02} \,, \quad
  \phiz := \phi_{01} - \phi_{02}
  \non
  \\
  & \theta := \theta_1 - \theta_2
  \aand
  \phi := \phi_1 - \phi_2 \,.
  \non
\Esist
Assuming that $\ell$ is a constant,
we write \eqref{prima} for both solutions
and integrate the difference with respect to time.
We obtain
\Beq
  \theta + \ell\phi + A^{2r} (1*\theta)
  = 1*f + \thetaz + \ell\phiz
  \quad \aeQ.
  \label{perdiffprima}
\Eeq
At this point, we multiply the above equality by $\theta$ and integrate over~$Q_t$, 
with an arbitrary $t\in(0,T)$.
On account of \eqref{propA} and~\eqref{propAB}, we~get
\begin{align}
  & \intQt |\theta|^2
  + \ell \intQt \phi\theta
  + \frac 12 \, \norma{A^r(1*\theta)(t)}^2
  \non
  \\
  & = \intQt (1*f) \theta
  + \intQt (\thetaz + \ell\phiz) \theta .
  \label{diffprima}
\end{align}
At the same time, we write \eqref{seconda} for both solution
and choose $v=\phi_2(t)$ and $v=\phi_1(t)$ in the inequalities we obtain, respectively.
Then, we sum up and integrate with respect to time.
As the contributions involving $\Beta$ cancel each other, we obtain
\begin{align}
  & \frac 12 \, \norma{\phi(t)}^2
  + \iot \norma{B^\sigma \phi(s)}^2 \, ds
  \non
  \\
  & \leq \frac 12 \, \norma\phiz^2
  - \intQt \bigl( \pi(\phi_1) - \pi(\phi_2) \bigr) \phi
  + \ell \intQt \theta\phi .
  \label{diffseconda}
\end{align}
Now, we add \eqref{diffprima} \pier{to \eqref{diffseconda}
and notice that the terms containing $\ell$ disappear}.
With the help of the \Lip\ continuity of~$\pi$
(see \eqref{hpPi}) and of the Young inequality,
we deduce that
\begin{align}
  & \intQt |\theta|^2
  + \frac 12 \, \norma{A^r(1*\theta)(t)}^2
  + \frac 12 \, \norma{\phi(t)}^2
  + \iot \norma{B^\sigma \phi(s)}^2 \, ds
  \non
  \\
  & \leq \frac 12 \, \intQt |\theta|^2
  + \intQt |1*f|^2
  + \intQt |\thetaz+\ell\phiz|^2
  + \frac 12 \, \norma\phiz^2
  + c \intQt |\phi|^2 .
  \label{percontdep}
\end{align}
Thus, \eqref{contdep} immediately follows by applying the Gronwall lemma.
\QED

In the nonlinear case of Theorem~\ref{Nonlinearell},
the equality \eqref{perdiffprima} has to be replaced~by
\begin{align}
  & \theta + \Ell(\phi_1) - \Ell(\phi_2) + A^{2r} (1*\theta)
  = 1*f + \thetaz + \Ell(\phi_{01}) - \Ell(\phi_{02})
  \non
  \\  
  & 
  \quad \hbox{where} \quad
  \gianni{\Ell(s) := \int_0^s \ell(\tau) \, d\tau
  \quad \hbox{for $s\in\erre$}.}
  \non
\end{align}
Moreover, the last term of \eqref{diffseconda} has to be modified in an obvious way.
Hence, the cancellation of the integrals involving $\ell$ does not occur any longer in summing up,
and the main difference with respect to the previous case is the following:
as a further contribution to the \rhs\ of the final inequality, we have the integral over $Q_t$ of the~sum
\begin{align}
  & \bigl( \pier{\Ell(\phi_2) - \Ell(\phi_1)} \bigr) \, \theta
  + \bigl( \theta_1 \, \ell(\phi_1) - \theta_2 \, \ell(\phi_2) \bigr) \phi
  \non
  \\
  & = \theta_1 \Bigl(
    \Ell(\phi_2) - \Ell(\phi_1) - \ell(\phi_1) (\phi_2-\phi_1)
  \Bigr)
  + \theta_2 \Bigl(
    \Ell(\phi_1) - \Ell(\phi_2) - \ell(\phi_2) (\phi_1-\phi_2)
  \Bigr).
  \non
\end{align}
However, this can be treated \pier{with the help of} our assumptions.
We write the Taylor expansion of $\Ell$ around any point $s\in\erre$
and see that \eqref{hpell} implies 
\Beq
  |\Ell(r) - \Ell(s) - \ell(s) (r-s)| \leq c \, |r-s|^2
  \quad \hbox{for every $r,s\in\erre$}.
  \non
\Eeq
Hence, we deduce that
\Beq
  \intQt \Bigl(
    \bigl(  \pier{\Ell(\phi_2) - \Ell(\phi_1)} \bigr) \, \theta
    + \bigl( \theta_1 \, \ell(\phi_1) - \theta_2 \, \ell(\phi_2) \bigr) \phi
  \Bigr)
  \leq c \intQt (|\theta_1| + |\theta_2) |\phi|^2 \,.
  \non
\Eeq
At this point, we invoke \eqref{hpsobolev} 
and apply~\eqref{compactineq}.
Thus, we can \pier{estimate the \rhs\ of} the above inequality and obtain
\begin{align}
  &\pier{\intQt \Bigl(
    \bigl(  \pier{\Ell(\phi_2) - \Ell(\phi_1)} \bigr) \, \theta
    + \bigl( \theta_1 \, \ell(\phi_1) - \theta_2 \, \ell(\phi_2) \bigr) \phi
  \Bigr)}
  \non
  \\
  & \leq c \iot \norma{|\theta_1(s)|+|\theta_2(s)|}_p \, \norma{\phi(s)}_q^2 \, ds
  \non
  \\
  & \leq c \, \bigl( \norma{\theta_1}_{\L\infty{\VA r}} + \norma{\theta_2}_{\L\infty{\VA r}} \bigr)
    \iot \norma{\phi(s)}_q^2 \, ds
  \non
  \\
  & \leq \frac 12 \iot \norma{B^\sigma\phi(s)}^2 \, ds 
  + c  \iot \norma{\phi(s)}^2 \, ds 
  \non
\end{align}
where the last value of $c$ also depends on the norms of $\theta_1$ and $\theta_2$ just written.
Therefore, we can come back to the modified \eqref{percontdep} 
and conclude as in the previous proof by applying the Gronwall lemma.

%%%%%%%%%%%%%%%%%%%%%%%%%%%%%%%%%%%%%%%%%%%%%%%%%%%%%%%%%%%%%%%%%%%%%%%%

\section{Approximation}
\label{APPROX}
\setcounter{equation}{0}

In this section, we prepare some auxiliary material that will be used
to perform the proofs of the existence parts of Theorems~\ref{Wellposedness} and~\ref{Nonlinearell} of the next section.
We introduce an approximating problem by fixing $\eps>0$
and replacing the function $\Beta$ and its subdifferential $\beta$
by their Moreau-Yosida regularizations $\Betaeps$ and $\betaeps$ at the level~$\eps$
(see, e.g., \cite[p.~28 and Prop.~2.11 p.~39]{Brezis}).
Thus, $\betaeps$~is monotone and \Lip\ continuous and coincides with the derivative of~$\Betaeps$.
Moreover, by also accounting for~\eqref{hpBeta}, it holds~that
\Beq
  0 \leq \Betaeps(r) \leq \Beta(r) , \quad
  \Betaeps(r) \leq \Beta_{\eps'}(r) \quad \hbox{if $\eps'<\eps$}, 
  \aand
  \lim_{\eps\searrow0} \Betaeps(r) = \Beta(r)
  \label{moreau}
\Eeq
for every $r\in\erre$.
We set for convenience
\Beq
  \Feps := \Betaeps + \Pi , \quad 
  \hbox{whence} \quad
  \Feps' = \betaeps + \pi .
  \label{defFeps}
\Eeq
Hence, the approximating problem consists in finding a pair $(\thetaeps,\phieps)$ satisfying
\begin{align}
  & \thetaeps \in \H1H \cap \L\infty{\VA r} \cap \L2{\VA{2r}}
  \qquad
  \label{regthetaeps}
  \\
  & \phieps \in \H1H \cap \L\infty{\VB\sigma} \cap \L2{\VB{2\sigma}}
  \label{regphieps}
\end{align}
\Accorpa\Regsoluzeps regthetaeps regphieps
and solving the system
\begin{align}
  & \dt\thetaeps + \ell(\phieps) \dt\phieps + A^{2r} \thetaeps = f
  \quad \aeQ
  \label{primaeps}
  \\[1mm]
  & \dt\phieps + B^{2\sigma}\phieps + \Feps'(\phieps) = \thetaeps \ell(\phieps)
  \quad \aeQ
  \label{secondaeps}
  \\[1mm]
  & \thetaeps(0) = \thetaz
  \aand
  \phieps(0) = \phiz \,.
  \label{cauchyeps}
\end{align}
\Accorpa\Pbleps primaeps cauchyeps
Notice that we have approximated the strong form \eqref{Iseconda} rather than~\eqref{seconda}.
The aim of this section is to prove that the above problem is well-posed.
We first treat the case that $\ell$ is a constant.
The \pier{more general} situation is considered later on in the section.

\Bthm
\label{Wellposednesseps}
Under the assumptions of Theorem~\ref{Wellposedness},
the approximating problem \Pbleps\ has a unique solution $(\thetaeps,\phieps)$
satisfying \Regsoluzeps.
\Ethm

As far as uniqueness is concerned, it suffices to observe that
\eqref{secondaeps} implies the analogue of \eqref{seconda}
obtained by replacing $\Beta$ by $\Betaeps$
and that $\Betaeps$ satisfies~\eqref{hpBeta}.
Thus, we can appy what we have just established in Section~\ref{UNIQUENESS}.
In order to prove the existence part, we use a Faedo-Galerkin scheme
depending on the parameter $n\in\enne$ and then we let $n$ tend to infinity.
By recalling \accorpa{eigen}{complete}, we introduce the subspaces
\Beq
  \VAn := \mathop{\rm span}\graffe{e_1,\dots,e_n}
  \aand
  \VBn := \mathop{\rm span}\graffe{\eta_1,\dots,\eta_n}
  \label{defVn}
\Eeq
and look for a pair $(\thetan,\phin)\in\H1{\VAn\times\VBn}$ satisfying
\begin{align}
  & \bigl( \dt\thetan + \ell \dt\phin + A^{2r} \thetan , v \bigr)
  = (f,v)
  \non
  \\
  & \quad \hbox{\aet\ for every $v\in\VAn$}
  \label{priman}
  \\[1mm]
  & \bigl( \dt\phin + B^{2\sigma}\phin + \Feps'(\phin) , v \bigr)
  = \ell (\thetan , v)
  \non
  \\
  & \quad \hbox{\aet\ for every $v\in\VBn$}
  \label{secondan}
  \\[1mm]
  & (\thetan(0),v) = (\thetaz,v)
  \aand
  (\phin(0),v) = (\phiz,v)
  \non
  \\
  & \quad \hbox{for every $v\in\VAn$ and $v\in\VBn$, respectively} \,.
  \label{cauchyn}
\end{align}
\Accorpa\Pbln priman cauchyn
Since $\eps$ is fixed at the moment, we did not stress 
the dependence of $(\thetan,\phin)$ on~$\eps$ in the notation.
First of all, we establish the existence of a global solution to the above problem.
We represent $(\thetan,\phin)$ in terms of the bases of $\VAn$ and~$\VBn$ as follows:
\Beq
  \thetan(t) = \somma j1n \thetan_j(t) e_j
  \aand
  \phin(t) = \somma j1n \phin_j(t) \eta_j
  \non
\Eeq
where the functions $\thetan_j$ and $\phin_j$ are looked for in $\H1\erre$.
If we equivalently \pier{let} $v=e_i$ in \eqref{priman} and $v=\eta_i$ in \eqref{secondan} with $i=1,\dots,n$,
we see that the system \accorpa{priman}{secondan} becomes
\begin{align}
  & \Bigl(
    \textstyle\somma j1n \dt\thetan_j \, e_j
    + \ell\textstyle\somma j1n \dt\phin_j \, \eta_j
    + \textstyle\somma j1n \lambda_j^{2r}\thetan_j e_j
  \ ,\ e_i
  \Bigr)
  = (f,e_i)
  \non
  \\
  & \quad \hbox{\aet\ for $i=1,\dots,n$}
  \non
  \\[1mm]
  & \Bigl(
    \textstyle\somma j1n \dt\phin_j \, \eta_j
    + \textstyle\somma j1n \mu_j^{2\sigma}\phin_j \eta_j
    + \Feps'(\textstyle\somma j1n \phin_j \eta_j)
    , \eta_i
  \Bigr)
  = \ell \Bigl( \textstyle\somma j1n \thetan_j e_j \ ,\ \eta_i \Bigr)
  \qquad
  \non
  \\
  & \quad \hbox{\aet\ for $i=1,\dots,n$}.
  \non
\end{align}
So, by introducing the $n$-column vectors
$\Theta:={}^t[\thetan_1,\dots,\thetan_n]$ and $\Phi:={}^t[\phin_1,\dots,\phin_n]$,
we obtain the following compact form of the system
\begin{align}
  & \Theta' + \calE \Phi' + \Lambda\Theta = g
  \aand
  \Phi' + M\Phi + \calF(\Theta,\Phi) = 0
  \qquad \hbox{or}
  \non
  \\
  & \Theta' - \calE \bigl( M\Phi + \calF(\Theta,\Phi) \bigr) + \Lambda\Theta = g
  \aand
  \Phi' + M\Phi + \calF(\Theta,\Phi) = 0
  \label{odesystem}
\end{align}
where the matrices $\calE$, $\Lambda$ and $M$ and the functions 
$g\in\L2{\erren}$ and $\calF:(\erren)^2\to\erren$
are defined~by
\begin{align}
  & \calE := \ell \, [(\eta_j,e_i)]_{i,j=1,\dots,n} \,, \quad
  \Lambda := \diag (\lambda_1^{2r},\dots,\lambda_n^{2r}) \,, \quad
  M := \diag (\mu_1^{2\sigma},\dots,\mu_n^{2\sigma})
  \qquad 
  \non
  \\[1mm]
  & g(t) := [(f(t),e_i)]_{i=1,\dots,n}
  \quad \aat
  \non
  \\
  & \calF(r,s)
  := \Biggl[
    \Biggl(
      \Feps' \Bigl(
        \textstyle\somma k1n s_k \eta_k
      \Bigr)
      - \ell \textstyle \somma j1n r_j e_j
      \ ,\ \eta_i
    \Biggr)
  \Biggr]_{i=1,\dots,n}
  \non
  \\
  & \quad \hbox{for $r=(r_1,\dots,r_n)\in\erren$ and $s=(s_1,\dots,s_n)\in\erren$} .
  \non
\end{align}
Since the \Lip\ continuity of~$\Feps'$ 
(see \eqref{defFeps})
implies the same property for~$\calF$
and the function $g$ belongs to $\L2\erren$,
every Cauchy problem for \eqref{odesystem} has a unique global solution
$(\Theta,\Phi)\in\H1{\erren\times\erren}$.
On the other hand, \eqref{cauchyn} yields an initial condition for $(\Theta,\Phi)$.
Therefore, coming back to problem \Pbln,
we conclude that it has a unique solution 
$(\thetan,\phin)\in\H1{\VAn\times\VBn}$.

At this point, we are ready to prove the existence part of Theorem~\ref{Wellposednesseps}.
This will be done by performing a number of a priori estimates
and \pier{passing to the limit by} compactness arguments.

\step
First a priori estimate

We test \eqref{priman} written at the time $s$ by $v=\thetan(s)$ 
and integrate with respect to $s$ over~$(0,t)$,
with an arbitrary $t\in(0,T)$.
In the same way, we test \eqref{secondan} by~$\dt\phin$, integrate with respect to time
and add the same quantity $\intQt\phin\dt\phin$ to both sides.
Then, we sum up and observe that the terms involving $\ell$ cancel each other.
Hence, on account of \accorpa{propA}{propAB}, we obtain
\begin{align}
  & \frac 12 \, \norma{\thetan(t)}^2
  + \iot \norma{A^r\thetan(s)}^2 \, ds
  + \iot \norma{\dt\phin(s)}^2 \, ds
  + \frac 12 \, \norma{\phin(t)}_{B,\sigma}^2
  + \iO \Betaeps(\phin(t))
  \qquad
  \non
  \\
  & = \frac 12 \, \norma{\thetan(0)}^2
  + \frac 12 \, \norma{B^\sigma\phin(0)}^2
  + \iO \Betaeps(\phin(0))
  \non
  \\
  & \quad {}
  + \iot (f(s),\thetan(s)) \, ds
  + \iot \bigl( \phin(s) - \pi(\phin(s)) , \dt\phin(s) \bigr) \, ds .
  \label{perprimastiman}
\end{align}
By also recalling \eqref{moreau}, we see that all the terms on the \lhs\ of \eqref{perprimastiman} are nonnegative.
The sum of the last two integrals on the \rhs\ is estimated\pier{,} owing to assumption \eqref{hpPi} and the Young inequality\pier{, as follows:} 
\begin{align}
  &\pier{\iot (f(s),\thetan(s)) \, ds
  + \iot \bigl( \phin(s) - \pi(\phin(s)) , \dt\phin(s) \bigr) \, ds}\non
\\
  &\leq \norma f_{\L2H}^2
  + \iot \norma{\thetan(s)}^2 \, ds
  + \frac 12 \iot \norma{\dt\phin(s)}^2 \, ds
  + c \iot \norma{\phin(s)}^2 \, ds
  + c \,.
  \non
\end{align}
Concerning the other terms, we observe that $\thetan(0)$ and $\phin(0)$ are the $H$-projections
of $\thetaz$ and $\phiz$ on $\VAn$ and~$\VBn$, respectively, due to~\eqref{cauchyn}.
By also accounting for the \Lip\ continuity of~$\betaeps$, we obtain for two of~them
\Beq
  \norma{\thetan(0)}^2
  + \iO \Betaeps(\phin(0))
  \leq \norma\thetaz^2 
  + c_\eps \bigl( \norma{\phin(0)}^2 + 1 \bigr)
  \leq \norma\thetaz^2 
  + c_\eps \bigl( \norma\phiz^2 + 1 \bigr)
  \leq c_\eps \,.
  \non
\Eeq
Finally, we have that
\Beq
  B^\sigma\phin(0)
  = B^\sigma \somma j1n (\phiz,\eta_j) \eta_j
  = \somma j1n \mu_j^\sigma (\phiz,\eta_j) \eta_j
  \non
\Eeq
whence also
\Beq
  \norma{B^\sigma\phin(0)}^2
  = \somma j1n |\mu_j^\sigma (\phiz,\eta_j)|^2
  \leq \somma j1\infty |\mu_j^\sigma (\phiz,\eta_j)|^2
  = \norma{B^\sigma\phiz}^2 \,.
  \non
\Eeq
Thus, coming back to \eqref{perprimastiman} and applying the Gronwall lemma, we conclude that
\Beq 
  \norma\thetan_{\L\infty H\cap\L2{\VA r}}
  + \norma\phin_{\H1H\cap\L\infty{\VB\sigma}}
  \leq c_\eps \,.
  \label{primastiman}
\Eeq

\step
Second a priori estimate

We test \eqref{priman} by $\dt\thetan$ and integrate with respect to time as before.
\pier{Thanks} to \accorpa{propA}{propAB} once more, we obtain
\Beq
  \iot \norma{\dt\thetan(s)}^2 \, ds
  + \frac 12 \, \norma{A^r\thetan(t)}^2
  = \frac 12 \, \norma{A^r\thetan(0)}^2
  + \iot \bigl( f(s) - \ell \dt\phin(s) , \dt\thetan(s) \bigr) \, ds .
  \non
\Eeq
By arguing as before in order to estimate the first term on the \rhs\
and using the Young inequality and \eqref{primastiman} for the second one,
we immediately conclude that
\Beq
  \norma\thetan_{\H1H\cap\L\infty{\VA r}}
  \leq c_\eps \,.
  \label{secondastiman}
\Eeq

\step
Limit

By \accorpa{primastiman}{secondastiman} and standard weak star compactness results,
we have for~a (not relabeled) subsequence
\Bsist
  & \thetan \to \thetaeps
  & \quad \hbox{weakly star in $\H1H\cap\L\infty{\VA r}$},
  \label{convthetan}
  \\
  & \phin \to \phieps
  & \quad \hbox{weakly star in $\H1H\cap\L\infty{\VB\sigma}$} .
  \label{convphin}
\Esist
\pier{In view of} the compact embeddings \eqref{compact} and applying a proper strong compactness result
(see, e.g., \cite[Sect.~8, Cor.~4]{Simon}),
we deduce that 
\Beq
  \thetan \to \thetaeps
  \aand
  \phin \to \phieps
  \quad \hbox{strongly in $\C0H$}.
  \label{strongn}
\Eeq
This implies, in particular, that $\Feps'(\phin)$ converges to $\Feps'(\phieps)$ in the same topology,
just by \Lip\ continuity.
We want to deduce that the following integrated version of the approximating problem 
\begin{align}
  & \ioT \bigl( \dt (\thetaeps + \ell\phieps)(s) - f(s), v(s) \bigr) \, ds
  + \ioT \bigl( A^r \thetaeps(s) , A^r v(s) \bigr) \, ds
  = 0,
  \label{intprimaeps}
  \\
  & \ioT \bigl( \dt\phieps(s) + \Feps'(\phieps(s)) - \ell \thetaeps(s) , v(s) \bigr) \, ds
  + \ioT \bigl( B^\sigma \thetaeps(s) , B^\sigma v(s) \bigr) \, ds
  = 0
  \label{intsecondaeps}
\end{align}
is fulfilled for every $v\in\L2{\VA r}$ and every $v\in\L2{\VB\sigma}$, respectively.
We start from the following integrated versions 
of the equations \eqref{priman} and~\eqref{secondan}\pier{:}
\begin{align}
  & \ioT \bigl( \dt (\thetan + \ell\phin)(s) - f(s), v(s) \bigr) \, ds
  + \ioT \bigl( A^r \thetan(s) , A^r v(s) \bigr) \, ds
  = 0 ,
  \label{intpriman}
  \\
  & \ioT \bigl( \dt\phin(s) + \Feps'(\phin(s)) - \ell \thetan(s) , v(s) \bigr) \, ds
  + \ioT \bigl( B^\sigma \thetan(s) , B^\sigma v(s) \bigr) \, ds
  = 0\pier{{},}
  \label{intsecondan}
\end{align}
which obviously hold for every $v\in\L2{\VAn}$ and every $v\in\L2{\VBn}$,
respectively, due to \accorpa{propA}{propB}.
We fix $m\in\enne$ and take any $v\in\L2\VAm$.
For every $n\geq m$, we have that $\VAm\subset\VAn$,
whence \eqref{intpriman} holds for~$v$.
By arguing similarly for \eqref{intsecondan} and then letting $n$ tend to infinity,
we deduce that \accorpa{intprimaeps}{intsecondaeps} are satisfied
for every $v\in\L2\VAm$ and $v\in\L2\VBm$, respectively.
By a simple density argument, we conclude that the same equations hold true 
for every $v\in\L2{\VA r}$ and every $v\in\L2{\VB\sigma}$, respectively.
We deduce that $(\thetaeps,\phieps)$ solves the equivalent problem
\begin{align}
  & \bigl( \dt (\thetaeps + \ell\phieps) - f, v \bigr) 
  + \bigl( A^r \thetaeps , A^r v \bigr)
  = 0
  \quad \hbox{\aet, for every $v\in\VA r$}\, , 
  \label{weakprimaeps}
  \\
  & \bigl( \dt\phieps + \Feps'(\phieps) - \ell \thetaeps , v \bigr) 
  + \bigl( B^\sigma \thetaeps , B^\sigma v \bigr)
  = 0
  \quad \hbox{\aet, for every $v\in\VB\sigma$} \, ,
  \label{weaksecondaeps}
\end{align}
so that \accorpa{primaeps}{secondaeps} follow from the lemma given below.
Finally, as \eqref{strongn} implies that $\thetan(0)$ and $\phin(0)$
converge to $\thetaeps(0)$ and $\phieps(0)$ strongly in~$H$,
we see that the initial conditions \eqref{cauchyeps} follow 
from the theory of orthogonal projections,
and the proof of Theorem~\ref{Wellposednesseps} is complete.
\QED

\Blem
\label{Strongsoluz}
Assume that $u\in\VA r$ and $\psi\in H$ satisfy
\Beq
  (A^r u , A^r v) = (\psi,v)
  \quad \hbox{for every $v\in\VA r$}.
  \label{hpstrongsoluz}
\Eeq
Then, we have that
\Beq
  u \in \VA{2r}
  \aand
  A^{2r}u = \psi.
  \label{strongsoluz}
\Eeq
Moreover, the same result hold if $A$ and $r$ 
are replaced by $B$ and~$\sigma$, respectively.
\Elem

\Bdim
We have that
\Beq
  u = \somma j1\infty (u,e_j) e_j \,, \quad
  \psi = \somma j1\infty (\psi,e_j) e_j 
  \aand
  \somma j1\infty |\lambda_j^r (u,e_j)|^2
  + \somma j1\infty |(\psi,e_j)|^2 < +\infty .
  \non
\Eeq
\pier{Moreover,} \eqref{hpstrongsoluz} written with $v=e_i$ implies that
$\lambda_i^{2r}(u,e_i)=(\psi,e_i)$ for every~$i$,
whence also $\somma j1\infty |\lambda_j^{2r} (u,e_j)|^2<+\infty$.
Thus, both conditions \eqref{strongsoluz} immediately follow.
\Edim

Now, we consider the case of a nonlinear function~$\ell$.
To this concern, we have the analogue of Theorem~\ref{Wellposednesseps}, namely

\Bthm
\label{Nonlinearelleps}
Under the assumptions of Theorem~\ref{Nonlinearell},
the approximating problem \Pbleps\ has a unique solution $(\thetaeps,\phieps)$
satisfying \Regsoluzeps.
\Ethm

\Bdim
The argument follows the same line of the previous proof 
and we just stress the modifications that are needed.
Concerning the choice of the Faedo-Galerkin scheme, it suffices to replace
$\ell$ by $\ell(\phin)$ in equations \eqref{priman} and~\eqref{secondan}.
However, the existence of a global discrete solution is no longer clear.
Indeed, the corresponding system of ordinary differential equations
still takes the form \eqref{odesystem}, 
but $\calE$ has to be replaced by $\calE(\Phi)$,
i.e., by a function depending on~$\Phi$, and the definition of $\calF$ has to be modified.
Precisely, we have to~set
\begin{align}
  & \calE(s) := \Biggl[ \Biggl( \ell \Bigl( \textstyle\somma j1n s_k \eta_k \Bigr) \eta_j , e_i \Biggr) \Biggr]_{i,j=1,\dots,n}
  \non
  \\
  & \calF(r,s) := \Biggl[ \Biggl(
    \Feps' \Bigl( \textstyle\somma k1n s_k \eta_k \Bigr)
    - \ell \Bigl( \textstyle\somma k1n s_k \eta_k \Bigr) \textstyle\somma j1n r_j e_j , \eta_i
  \Biggr) \Biggr]_{i=1,\dots,n}
  \non
  \\
  & \quad \hbox{for $r=(r_1,\dots,r_n)\in\erren$ and $s=(s_1,\dots,s_n)\in\erren$} .
  \non
\end{align}
Thus, if $\ell'$ does not vanish identically,
the derivatives of $\calF$ with respect to $s_j$ have a linear growth with respect to~$r$
and the \Lip\ condition that is needed to ensure the existence of a global solution fails.
For this reason, we can only conclude that 
the system has a unique maximal solution $(\Theta,\Phi)$ defined in some interval $[0,T_n)\subseteq[0,T]$.
Thus, in principle, the maximal solution $(\thetan,\phin)$ to the discrete problem \Pbln\ 
exists and is defined in the same interval.
However, the estimates we can perform show that, for every~$n$, 
$(\thetan,\phin)$ is bounded in $\L\infty{\VAn\times\VBn}$
by a constant that depends on $T$ but not on~$T_n$.
Thus, the same happens for $(\Theta,\Phi)\in\L\infty{\erren\times\erren}$,
so that the general \pier{theory} of ordinary differential equations ensures \pier{the existence of a solution 
in the whole interval $[0,T]$}.

Coming to the estimates, the first one goes exactly as we did to obtain~\eqref{primastiman}.
Indeed, by testing the equations in the same way, 
the terms involving $\ell$ cancel each other also in the nonlinear case.
Concerning the second estimate, we only have to recall that $\ell$ is bounded,
so that $\ell(\phin)\dt\phin$ is bounded in $\L2H$ exactly as $\ell\dt\phin$ was before.
So, the analogue of \eqref{secondastiman} follows.
Furthermore, in taking the limit as $n$ tends to infinity, 
it suffices to add to the above argument the following convergence property:
$\ell(\phin)$ converges to $\ell(\phieps)$ strongly in $\C0H$,
as a consequence of \eqref{strongn} and the \Lip\ continuity of~$\ell$.
\pier{Recalling \eqref{convphin} as well, it is easy to see} that $\ell(\phin)\dt\phin$ and $\ell(\phin)\thetan$ 
\pier{are bounded \gianni{in $\L2H$} and}
converge weakly in \gianni{($\L2\Luno$\piernew{,} thus~in)} $\L2H$ to $\ell(\phieps)\dt\phieps$ and $\ell(\phieps)\thetaeps$, respectively,
so that we can take the limit of the corresponding terms 
in the analogues of \eqref{intpriman} and~\eqref{intsecondan}.
Finally, Lemma~\ref{Strongsoluz} applies also in the present case.
This completes the proof.
\Edim

%%%%%%%%%%%%%%%%%%%%%%%%%%%%%%%%%%%%%%%%%%%%%%%%%%%%%%%%%%%%%%%%%%%%%%%%

\section{Existence and regularity}
\label{EXISTENCE}
\setcounter{equation}{0}

In this \pier{section}, we conclude the proofs of the \pier{existence and regularity results
stated} in Section~\ref{STATEMENT}.
In order to establish the existence parts of both Theorems~\ref{Wellposedness} and~\ref{Nonlinearell},
we assume the general hypothesis~\eqref{hpell} and avoid using~\eqref{hpsobolev}.
We start from the approximating problem \Pbleps\
(whose unique solution exists thanks to Theorems~\ref{Wellposednesseps} and~\ref{Nonlinearelleps})
and perform some a~priori estimates.

\step
\pier{Uniform estimates}

We multiply \eqref{primaeps} and \eqref{secondaeps} by $\thetaeps$ and~$\dt\phieps$, respectively,
sum up and integrate over~$Q_t$, with an arbitrary $t\in(0,T)$.
We notice that the terms involving $\ell$ cancel each other 
and add the same term $\intQt\phieps\dt\phieps$ to both sides.
We obtain
\begin{align}
  & \frac 12 \, \norma{\thetaeps(t)}^2
  + \iot \norma{A^r\thetaeps(s)}^2 \, ds
  + \intQt |\dt\phieps|^2
  + \frac 12 \, \norma{\phieps(t)}_{B,\sigma}^2
  + \iO \Betaeps(\phieps(t))
  \non
  \\
  & = \frac 12 \, \norma\thetaz^2
  + \frac 12 \, \norma\phiz_{B,\sigma}^2
  + \iO \Betaeps(\phiz)
  + \intQt \bigl( f + \phieps - \pi(\phieps) \bigr) \dt\phieps \,.
  \label{perprimastima}
\end{align}
We treat the \rhs\
by first using the Young inequality and \eqref{moreau}\pier{, then invoking}
the assumptions \HPdati\ and the linear growth of~$\pi$.
\pier{Hence,} by applying the Gronwall lemma, we \pier{infer} that
\Beq
  \norma\thetaeps_{\L\infty H\cap\L2{\VA r}}
  + \norma\phieps_{\H1H\cap\L\infty{\VB\sigma}}
  + \norma{\Betaeps(\phieps)}_{\L\infty\Luno}
  \leq c \,.
  \label{primastima}
\Eeq
\pier{Next, we} multiply \eqref{primaeps} by $\dt\thetaeps$ and integrate over~$Q_t$ as before.
We obtain
\Beq
  \intQt |\dt\thetaeps|^2
  + \frac 12 \, \norma{A^r\thetaeps(t)}^2
  = \pier{{}\frac 12 \, \norma{A^r\thetaz}^2 {}} + \intQt \bigl( f-\ell(\phieps)\dt\phieps \bigr) \dt\thetaeps \,.
  \non 
\Eeq
With the help of the boundedness of~$\ell$,
the Young inequality and~\eqref{primastima}, we conclude that
\Beq
  \norma\thetaeps_{\H1H\cap\L\infty{\VA r}}
  \leq c \,.
  \label{secondastima}
\Eeq

\step
Limit

As in the previous section, we owe to \pier{weak star} and strong compactness results
and deduce from \accorpa{primastima}{secondastima} that (at least for a subsequence)
\Bsist
  & \thetaeps \to \theta
  & \quad \hbox{weakly star in $\H1H\cap\L\infty{\VA r}$}
  \non
  \\
  && \quad \hbox{and strongly in $\C0H$}\pier{,}
  \label{convtheta}
  \\[2mm]
  & \phieps \to \phi
  & \quad \hbox{weakly star in $\H1H\cap\L\infty{\VB\sigma}$}
  \non
  \\
  && \quad \hbox{and strongly in $\C0H$}.
  \label{convphi}
\Esist
We deduce that the initial conditions \eqref{cauchy} hold true
and that $\pi(\phieps)\to\pi(\phi)$ and $\ell(\phieps)\to\ell(\phi)$
strongly in $\C0H$ by \Lip\ continuity.
It follows that \pier{$\ell(\phieps)\dt \phieps\to\ell(\phi)\dt\phi$ and}
$\ell(\phieps)\thetaeps\to\ell(\phi)\theta$ weakly in $\L2H$
\gianni{(by~boundedness in $\L2H$ and weak convergence in $\L2\Luno$)}.
At~this point, it is \sfw\ to show that $(\theta,\phi)$ satisfies
the integrated version of \eqref{prima} similar to~\eqref{intprimaeps}.
Then, it also satisfies the \pier{analog} of \eqref{weakprimaeps},
so that both the full regularity \eqref{regtheta} and the equation \eqref{prima} follow
from Lemma~\ref{Strongsoluz}.
So, it remains to prove that $\phi$ also satisfies \eqref{regBetaphi}
and that $(\theta,\phi)$ solves the variational inequality~\eqref{seconda}.
We recall that all the above convergence properties 
hold for some subsequence $\eps_n\searrow0$,
which we can assume to be strictly decreasing without loss of generality.
Moreover, we can assume that $\phi_{\eps_n}$ converges to $\phi$ \aeQ.
We first prove that
\Beq
  \intQ \Beta(\phi) 
  \leq \liminf_{n\to\infty} \intQ \Beta_{\eps_n}(\phi_{\eps_n}).
  \label{liminfBeta}
\Eeq
We take arbitrary indices $n$ and $m$ with $n>m$.
Then, $\eps_n<\eps_m$ and we can apply~\eqref{moreau}.
We deduce that
\Beq
  \Beta_{\eps_m}(\phi_{\eps_n})
  \leq \Beta_{\eps_n}(\phi_{\eps_n})
  \quad \hbox{\aeQ, for every $n>m$}.
  \non
\Eeq
Since $\Beta_{\eps_m}$ is (\Lip) continuous, we thus have that
\Beq
  \Beta_{\eps_m}(\phi)
  = \lim_{n\to\infty} \Beta_{\eps_m}(\phi_{\eps_n})
  = \liminf_{n\to\infty} \Beta_{\eps_m}(\phi_{\eps_n})
  \leq \liminf_{n\to\infty} \Beta_{\eps_n}(\phi_{\eps_n})
  \quad \aeQ.
  \non
\Eeq
On the other hand, the last \pier{condition in} \eqref{moreau} implies that
\Beq
  \Beta(\phi)
  = \lim_{m\to\infty} \Beta_{\eps_m}(\phi)
  \quad \aeQ.
  \non
\Eeq
Therefore, \eqref{liminfBeta} follows from the Fatou lemma.
At this point, we can easily conclude.
From one side, \eqref{liminfBeta} implies \eqref{regBetaphi} due to~\eqref{primastima}.
On the other hand, it is easily seen that \eqref{secondaeps}
implies the analogue of~\eqref{seconda}, i.e.,
the variational inequality obtained by replacing $\Beta$, $\theta$ and $\phi$ in \eqref{seconda}
by $\Betaeps$, $\thetaeps$ and~$\phieps$, respectively.
Therefore, by also combining the strong and weak convergence properties \accorpa{convtheta}{convphi} 
and writing $\eps$ instead of~$\eps_n$ for simplicity,
we have for every $v\in\L2{\VB\sigma}$
\begin{align}
  & \intQ \Beta(\phi)
  + \ioT \bigl( B^\sigma\phi(t) , B^\sigma (\phi(t)-v(t)) \bigr) \, dt
  \non
  \\
  & \leq \liminf_{\eps\searrow0} \intQ \Betaeps(\phieps)
  + \liminf_{\eps\searrow0} \ioT \bigl( B^\sigma\phieps(t) , B^\sigma (\phieps(t)-v(t)) \bigr) \, dt
  \non
  \\
  & \leq \liminf_{\eps\searrow0} \Bigl(
    \intQ \Betaeps(\phieps)
    + \ioT \bigl( B^\sigma\phieps(t) , B^\sigma (\phieps(t)-v(t)) \bigr) \, dt
  \Bigr)
  \non
  \\
  & \leq \liminf_{\eps\searrow0} \Bigl(
    \intQ \bigl( \ell(\phieps)\thetaeps - \dt\phieps - \pi(\phieps) \bigr) (\phieps-v)
    + \intQ \Betaeps(v)
  \Bigr)
  \non
  \\
  & = \intQ \bigl( \ell(\phi)\theta - \dt\phi - \pi(\phi) \bigr) (\phi-v)
  + \intQ \Beta(v)
  \non
\end{align}
In the last equality we have used the last \eqref{moreau} as well.
Hence, \eqref{seconda} is proved and the proof is complete.
\QED

\step
Proofs of Propositions \ref{Secondaqo} and \ref{Nonlinearsecondaqo}

Since we use just the boundedness of~$\ell$, the same proof holds for both propositions.
We start from the approximating problem once more
and, \aat, we write \eqref{secondaeps} at the time~$t$ in the form
\Beq
  \bigl( B^\sigma\phieps(t) , B^\sigma v \bigr) 
  + \iO \betaeps(\phieps(t)) \, v
  = \bigl( \thetaeps(t) \, \ell(\phieps(t)) - \dt\phi(t) - \pi(\phi(t)) , v \bigr)
  \quad \hbox{for every $v\in\VB\sigma$}.
  \non
\Eeq
\pier{Thanks} to \eqref{hpqo}, we can choose $v=\betaeps(\phieps(t))$ and obtain
\Beq
  \iO |\betaeps(\phi(t))|^2
  \leq \bigl( \thetaeps(t) \, \ell(\phieps(t)) - \dt\phieps(t) - \pi(\phieps(t)) , \betaeps(\pier{\phieps(t)}) \bigr) .
  \non
\Eeq
Then, the Young inequality, the boundedness of~$\ell$ and \eqref{primastima} immediately yield
\Beq
  \norma{\betaeps(\phieps)}_{\L2H}
  \leq c \,.
  \non
\Eeq
As a consequence, we also have from \eqref{secondaeps} that
\Beq
  \norma{B^{2\sigma}\phieps}_{\L2H}
  \leq c \,.
  \non
\Eeq
Therefore, coming back the the proof of Theorems~\ref{Wellposedness} and~\ref{Nonlinearell} just performed,
we~see that we can add to \accorpa{convtheta}{convphi} the further convergence properties
\Beq
  \betaeps(\phieps) \to \xi
  \aand
  B^{2\sigma}\phieps \to B^{2\sigma}\phi
  \quad \hbox{weakly in $\L2H$}.
  \label{convxi}
\Eeq
At this point, it is clear that $\phi\in\L2{\VB{2\sigma}}$ and that \eqref{secondaqo} holds true.
In order to \pier{check the inclusion property in}~\eqref{regxi}, 
it suffices to recall the strong convergence \eqref{convphi} of $\phieps$
and \pier{use the maximal monotonicity of $\beta$ by applying,} e.g., \cite[Lemma~2.3, p.~38]{Barbu}.
\QED

Finally, we prove Theorems~\ref{Regularity} and~\ref{Nonlinearregularity}.
We start from the first of them, i.e., we assume that $\ell$ is a constant.
Clearly, it suffices to establish the estimates 
corresponding to \gianni{\eqref{regdtphi}}
on the solution $(\thetaeps,\phieps)$ to the approximating problem \Pbleps.

\step
\pier{Regularity estimate}

We proceed formally, for brevity.
We differentiate \eqref{secondaeps} with respect to time
and test the \pier{resulting equality} by~$\dt\phieps$.
Then, we integrate over~$Q_t$, as usual.
We obtain
\begin{align}
  & \frac 12 \, \norma{\dt\phieps(t)}^2
  + \iot \norma{B^\sigma\dt\phieps(s)}^2
  + \intQt \betaeps'(\phieps) |\dt\phieps|^2
  \non
  \\
  & = \frac 12 \, \norma{\dt\phieps(0)}^2
  + \intQt \bigl( \ell\dt\thetaeps - \pi'(\phieps) \dt\phieps \bigr) \dt\phieps \,.
  \label{perprimareg}
\end{align}
The last integral can be trivially estimated \pier{owing to} 
the \Lip\ continuity of~$\pi$, the Young inequality and \accorpa{primastima}{secondastima} as follows
\Beq
  \intQt \bigl( \ell\dt\thetaeps - \pi'(\phieps) \dt\phieps \bigr) \dt\phieps 
  \leq c \intQt |\dt\phieps|^2
  + c \intQt |\dt\thetaeps|^2
  \leq c \,.
  \non
\Eeq
On the other hand, we have from \eqref{secondaeps} 
\Beq
  \norma{\dt\phieps(0)}
  \leq \ell \, \norma\thetaz
  + \norma{B^{2\sigma}\phiz}
  + \norma{\betaeps(\phiz)}
  + \norma{\pi(\phiz)}
  \leq \norma{\betaz(\phiz)} + c
  = c
  \non
\Eeq
since $|\betaeps(r)|\leq|\betaz(r)|$ for every $r\in D(\beta)$
(see, e.g., \cite[Prop.~2.6, p.~28]{Brezis}).
Hence, we conclude that
\Beq
  \norma{\dt\phieps}_{\L\infty H\cap\L2{\VB\sigma}}
  \leq c \,.
  \label{primareg}
\Eeq
\QED

\step
Proof of Theorem~\ref{Nonlinearregularity}

\gianni{We show how to modify the derivation of estimate~\eqref{primareg}
in the nonlinear case \eqref{hpell} by also assuming~\eqref{hpsobolev}.
%Concerning~\eqref{primareg}, we just replace $\ell$ by $\ell(\phieps)$ and use the boundedness of~$\ell$.
%On the contrary, \eqref{secondareg}~is more delicate.
The difference is the \rhs\ of \eqref{perprimareg} one obtains 
by differentiating \eqref{secondaeps} with respect to time and then testing by~$\dt\phieps$}.
Namely, it contains the more complicated terms
\Beq
  \intQt \dt\thetaeps \, \ell(\phieps) \dt\phieps
  + \intQt \thetaeps \, \ell'(\phieps) |\dt\phieps|^2 \,.
  \non
\Eeq
The first one is treated as before by using the boundedness of~$\ell$.
\pier{About the second one, we recall that $\ell'$ is bounded by~\eqref{hpell}; then,}
owing to the \Holder\ inequality and \eqref{hpsobolev},
and accounting for \eqref{secondastima} and the compacness inequality \eqref{compactineq}\pier{, we
have that}
\begin{align}
  & \intQt \thetaeps \, \ell'(\phieps) |\dt\phieps |^2
  \leq c \iot \norma{\thetaeps(s)}_p \, \norma{\dt\phieps(s)}_q^2 \, ds
  \leq c \iot \norma{\thetaeps(s)}_{A,r} \, \norma{\dt\phieps(s)}_q^2 \, ds
  \non
  \\
  & \leq c \iot \norma{\dt\phieps(s)}_q^2 \, ds
  \leq \frac 12 \iot \norma{B^\sigma\dt\phieps(s)}^2 \, ds
  + c \iot \norma{\dt\phieps(s)}^2 \, ds \,.
  \non
\end{align}
By combining this with the modified~\eqref{perprimareg} and applying the Gronwall lemma,
we obtain~\eqref{primareg}.
% Finally, estimates \eqref{secondareg} and \eqref{terzareg} can be proved as before,
% just by replacing $\ell$ by $\ell(\phieps)$ 
% and using the boundedness of $\ell$ postulated in~\eqref{hpell}.
Thus, the proof is complete.

%%%%%%%%%%%%%%%%%%%%%%%%%%%%%%%%%%%%%%%%%%%%%%%%%%%%%%%%%%%%%%%%%%%%%%%%

\section{Longtime \bhv}
\label{LONGTIME}
\setcounter{equation}{0}

In this section, we prove Theorem~\ref{Longtime}.
Thus, in addition to the assumption \HPstruttura\ on the structure
and \eqref{hpz} on the initial data,
we suppose that the other hypotheses of the statement are in force
and study the \omegalimit\ of the unique global solution.
In order to show that $\omega$ is \pier{nonempty} and for a further use, 
we need some global estimates on $(\theta,\phi)$ 
on the half line~$[0,+\infty)$.
In developping our argument, we start from the approximating problem once more,
which clearly has a unique solution $(\thetaeps,\phieps)$ defined in the whole of~$[0,+\infty)$.
We notice that the convergence of the approximating solution to $(\theta,\phi)$
we have proved in Section~\ref{EXISTENCE} 
holds for every~$T$ 
(and for the whole family, i.e.,~not only for a subsequence),
thank to our uniqueness result proved in Section~\ref{UNIQUENESS}.

\step
First global estimate

As we did to prove \eqref{primastima}, 
we test \eqref{primaeps} and \eqref{secondaeps} by $\thetaeps$ and~$\dt\phieps$, respectively,
integrate over $Q_t$ and sum up.
However, in contrast to the previous argument,
we avoid adding the same term to both sides of the equality we get.
Since the terms involving $\ell$ cancel each other and \eqref{hpcoerc} holds,
we obtain
\begin{align}
  & \frac 12 \, \norma{\thetaeps(t)}^2
  + \iot \norma{A^r\thetaeps(s)}^2 \, ds
  + \intQt |\dt\phieps|^2
  + \frac 12 \, \norma{B^\sigma\phieps(t)}^2
  + \alpha \, \norma{\phieps(t)}^2 - C
  \non
  \\
  & \leq \frac 12 \, \norma\thetaz^2
  + \frac 12 \, \norma{B^\sigma\phiz}^2
  + \iO \Betaeps(\phiz) 
  + \intQt f \thetaeps \,.
  \non
\end{align}
The first three terms on the \rhs\ are bounded uniformly with respect to~$\eps$
by the first \pier{condition in}~\eqref{moreau} and \eqref{hpz}.
Concerning the last one, we owe to \eqref{hpfbis} and have~that
\Beq
  \intQt f \thetaeps
  \leq \sup_{0\leq s\leq t} \norma{\thetaeps(s)} \iot \norma{f(s)} \, ds
  \leq \frac 14 \sup_{0\leq s\leq t} \norma{\thetaeps(s)}^2
  + \norma f_{\LL1H}^2 \,.
  \non
\Eeq 
At this point, it is \sfw\ to conclude that
\Beq
  \norma\thetaeps_{\LL\infty H}
  + \norma{A^r\thetaeps}_{\LL2H}
  + \norma\phieps_{\LL\infty{\VB\sigma}}
  + \norma{\dt\phieps}_{\LL2H}
  \leq c \,.
  \label{primaglobale}
\Eeq

\step
Second global estimate

We test \eqref{primaeps} by $\dt\thetaeps$.
By integrating over $Q_t$ and 
using the Young inequality, the boundedness of~$\ell$ and assumption \eqref{hpfbis} on~$f$, 
we obtain
\begin{align}
  & \intQt |\dt\thetaeps|^2
  + \frac 12 \, \norma{A^r\thetaeps(t)}^2
  = \intQt \bigl( f - \ell(\phieps) \dt\phieps \bigr) \dt\thetaeps 
  \non
  \\
  & \leq \frac 12 \intQt |\dt\thetaeps|^2
  + c \, \norma f_{\LL2H}^2 
  + c \, \intQt |\dt\phieps|^2 \,.
  \non
\end{align}
Hence, \eqref{primastima} immediately yields
\Beq
  \norma\thetaeps_{\LL\infty{\VA r}}
  + \norma{\dt\thetaeps}_{\LL2H}
  \leq c \,.
  \label{secondaglobale}
\Eeq

\step
Basic global estimates

By letting $\eps$ tend to zero in \accorpa{primaglobale}{secondaglobale}, 
we see that the solution $(\theta,\phi)$ we are studying enjoys the properties
\begin{align}
  & \theta \in \pier{\LL\infty {V_A^r}}
  \aand
  \phi \in \LL\infty{\VB\sigma} \subset \LL\infty H
  \label{globalbdd}
  \\
  & \int_0^{+\infty} \norma{A^r\theta(s)}^2 \, ds
  + \int_0^{+\infty} \norma{\dt\theta(s)}^2 \, ds
  + \int_0^{+\infty} \norma{\dt\phi(s)}^2 \, ds 
  < +\infty \,.
  \label{peromegalimit}
\end{align}
In particular, the \omegalimit\ $\omega$ is \pier{nonempty}.

The next step consists in proving the properties of the elements of $\omega$ we have stated.
Thus, we fix $(\thetao,\phio)\in\omega$ 
and a corresponding sequence $\{\tn\}$ as in the definition \eqref{defomegalimit},
and, for every~$n$, we study the limits on a fixed time interval~$(0,T)$
of the functions $\thetan$ and $\phin$ defined~by
\Beq
  \thetan(t) := \theta(t+\tn)
  \aand
  \phin(t) := \phi(t+\tn)
  \quad \hbox{for $t\in[0,T]$}.
  \label{defthetanphin}
\Eeq
The global estimates \accorpa{globalbdd}{peromegalimit} immediately yield~that
\begin{align}
  & \norma\thetan_{\L\infty H}
  + \norma\phin_{\L\infty{\VB\sigma}}
  \leq c
  \non
  \\
  & \lim_{n\to\infty} \Biggl(
    \ioT \norma{A^r\thetan(s)}^2 \, ds
    + \ioT \norma{\dt\thetan(s)}^2 \, ds
    + \ioT \norma{\dt\phin(s)}^2 \, ds     
  \Biggr)
  = 0 \,.
  \non
\end{align}
By standard compactness results it follows that
\begin{align}
  & \thetan \to \thetai
  \quad \hbox{weakly \pier{star} in $\H1H\cap\pier{\L\infty{\VA r}}$}
  \label{convthetani}
  \\
  & \phin \to \phii
  \quad \hbox{weakly star in $\H1H\cap\L\infty{\VB\sigma}$}
  \label{convphini}
\end{align}
at least for a subsequence, where $\thetai$ and $\phii$ satisfy
\Beq
  A^r \thetai = 0 , \quad
  \dt\thetai = 0
  \aand
  \dt\phii = 0
  \quad \aeQ .
  \label{opnulli}
\Eeq
In particular, both $\thetai$ and $\phii$ are time-independent,
so that we can define the elements $\thetas\in\VA r$ and $\phis\in\VB\sigma$
by setting
\Beq
  \thetas := \thetai(t)
  \aand
  \phis := \phii(t)
  \quad \hbox{for every $t\in[0,T]$}.
  \non
\Eeq
Our aim is to show that $(\thetao,\phio)=(\thetas,\phis)$
and that $(\thetas,\phis)$ solves \accorpa{primas}{secondas}.
By the weak convergence in $\C0H$ implied by the weak convergence in $\H1H$, we have that
\Beq
  \thetas
  = \thetai(0)
  = \lim_{n\to\infty} \thetan(0)
  = \lim_{n\to\infty} \theta(\tn)
  = \thetao
  \non
\Eeq
where the limits are understood in the weak topology of~$H$.
Similarly, we obtain that $\phis=\phio$.
As far as the limiting system is concerned,
equation \eqref{primas} follows from the first \pier{equality in}~\eqref{opnulli}.
It remains to prove~\eqref{secondas}.
To this end, we observe that the pair $(\thetan,\phin)$ obviously satisfies~\eqref{seconda},
whence also~\eqref{intseconda}.
On the other hand, we can invoke the compact embedding $\VB\sigma\subset H$ (see~\eqref{compact})
and apply, e.g., \cite[Sect.~8, Cor.~4]{Simon} to obtain
\Beq
  \phin \to \phii
  \quad \pier{\hbox{strongly in $\C0H$}.}
  \label{strongphini}
\Eeq
By \Lip\ continuity, it follows that
\Beq
  \pi(\phin) \to \pi(\phii)
  \quad \hbox{strongly in $\C0H$}.
  \label{convpi}
\Eeq
Now, we prove that
\Beq
  \thetan \, \ell(\phin) \to \thetai \, \ell(\phii)
  \quad \hbox{weakly \pier{star in $\L\infty H$}}
  \label{convprod}
\Eeq
by distinguishing the cases of the statement that concern~$\ell$.
If $\ell$ is a constant, then \eqref{convprod} trivially follows from~\eqref{convthetani}.
In the opposite case, we invoke \eqref{hpsobolevbis} and apply \cite[Sect.~8, Cor.~4]{Simon}.
Hence, from \accorpa{convthetani}{convphini} we deduce that
\begin{align}
  & \thetan \to \thetai
  \quad \hbox{weakly \pier{star in $\L\infty{\Lx p}$},} 
  \non
  \\
  & \phin \to \phii
  \quad \hbox{strongly in $\C0{\Lx q}$}, \quad \hbox{whence}
  \non
  \\
  & \ell(\phin) \to \ell(\phii)
  \quad \hbox{strongly in $\C0{\Lx q}$}.
  \non
\end{align}
Thus, \eqref{convprod} follows also in this case since $(1/p)+(1/q)=1/2$.
We remark that \eqref{convpi}, \eqref{convprod} and \eqref{strongphini} imply (at~least)
\Beq
  \pi(\phin) \, \phin \to \phi(\phii) \, \phii
  \aand
  \thetan \, \ell(\phin) \, \phin \to \thetai \, \ell(\phii) \, \phii
  \quad \hbox{weakly in $\LQ1$}.
  \non  
\Eeq
At this point,
we can easily let $n$ tend to infinity in \eqref{intseconda} written for $(\thetan,\phin)$.
By also accounting for the \pier{lower semicontinuity of the convex function $v\mapsto \intQ \Beta(v) $ 
in $L^2(Q)$,} we have~that
\begin{align}
  & \intQ \Beta(\phii)
  + \ioT \bigl( B^\sigma\phii(t) , B^\sigma (\phii(t)-v(t)) \bigr) \, dt
  \non
  \\
  \separa
  & \leq \liminf_{\eps\searrow0} \intQ \Beta(\phin)
  + \liminf_{n\to\infty} \ioT \bigl( B^\sigma\phin(t) , B^\sigma (\phin(t)-v(t)) \bigr) \, dt
  \non
  \\
  \separa
  & \leq \liminf_{n\to\infty} \Bigl(
    \intQ \Beta(\phin)
    + \ioT \bigl( B^\sigma\phin(t) , B^\sigma (\phin(t)-v(t)) \bigr) \, dt
  \Bigr)
  \non
  \\
  \separa
  & \leq \liminf_{n\to\infty} \Bigl(
    \intQ \bigl( \thetan\ell(\phin) - \dt\phin - \pi(\phin) \bigr) (\phin-v)
    + \intQ \Betaeps(v)
  \Bigr)
  \non
  \\
  \separa
  & = \intQ \bigl( \thetai\ell(\phii) - \dt\phii - \pi(\phii) \bigr) (\phii-v)
  + \intQ \Beta(v)
  \non
\end{align}
where we have kept $\dt\phii$ for clarity, even though it vanishes.
Thus, \eqref{intseconda} holds true for $(\thetai,\phii)$.
This implies that \eqref{seconda} holds as well.
But the latter coincides with~\eqref{secondas}
and the proof is complete.

%%%%%%%%%%%%%%%%%%%%%%%%%%%%%%%%%%%%%%%%%%%%%%%%%%%%%%%%%%%%%%%%%%%%%%%%

\section{Convergence to a phase relaxation problem}
\label{RELAXATION}
\setcounter{equation}{0}

In this section we discuss the asymptotic \bhv\ of the solution to our problem 
as $\sigma\searrow 0 $. 

We assume that \accorpa{hpconst}{defbetapi} are satisfied,
\begin{align}
\label{a-pier0}
&\hbox{$\ell \, $ is a constant},
\\
\label{a-pier1}
&\pi (v) = -\, \gamma \, v  \quad \hbox{for all } v\in \erre , \hbox{ with a fixed constant }  \gamma \geq 0,
\end{align}
and there exist some $\sigma_0 >0$ and some family of data 
$\{ f_{\sigma},\, \theta_{0,\sigma}, \, \phi_{0,\sigma} \}$ 
such that
\begin{align}
\label{a-pier2}
 &f_{\sigma} \to f  \quad \hbox{in } \, \L2H \quad \hbox{as }\, \sigma\searrow 0, \\
\label{a-pier3}
&\theta_{0,\sigma} \to \theta_0  \quad \hbox{in } \, \VA r \quad \hbox{as }\, \sigma\searrow   0,\\[1mm]
 &\phi_{0,\sigma}  \in \VB{\sigma} \quad \hbox{and} 
 \quad  \norma{\phi_{0,\sigma}}_{\VB{\sigma}} + \norma{\Beta(\phi_{0,\sigma})}_{L^1(\Omega)} 
 \leq c   
 \quad \hbox{for all } \sigma \in (0,\sigma_0], \non \\ 
 &\quad  \phi_{0,\sigma} \to \phi_0  \quad \hbox{in } \, H \quad \hbox{as }\, \sigma\searrow   0.\label{a-pier4}
\end {align}
Concerning \eqref{a-pier4}, we just note that if $\phiz \in \VB {\sigma_0} $ with 
$ \Beta(\phiz) \in \Luno $ (cf.~\eqref{hpz}), then the constant sequence $ \phi_{0,\sigma}= \phi_0 $ directly works in \eqref{a-pier4}.

We are dealing with the solution $(\theta_\sigma,\phi_\sigma)$ to the system (cf.~\eqref{prima}--\eqref{cauchy})
\begin{align}
  & \dt\theta_\sigma + \ell \, \dt\phi_\sigma + A^{2r} \theta_\sigma = f_\sigma
  \quad \aeQ \, ,
  \label{primasig}
  \\[1mm]
  & \bigl( \dt\phi_\sigma(t) , \phi_\sigma(t) - v \bigr)
  + \bigl( B^\sigma \phi_\sigma(t) , B^\sigma (\phi_\sigma(t) - v) \bigr) 
%  \non
%  \\
%  & \quad {}
  + \iO \Beta(\phi_\sigma(t))
  -\bigl( \gamma \, \phi_\sigma(t) , \phi_\sigma(t) - v \bigr)
  \non
  \\
  & \leq \bigl( \ell\, \theta_\sigma(t) , \phi_\sigma(t) - v \bigr)
  + \iO \Beta(v)
%  \non
%  \\
%  & 
  \quad 
  \hbox{\aat\ and every $v\in\VB\sigma$} \, ,
  \label{secondasig}
  \\[1mm]
  & \theta_\sigma(0) = \theta_{0,\sigma}
  \aand
  \phi_\sigma (0) = \phi_{0,\sigma} \,,
  \label{cauchysig}
\end{align}
where $(\theta_\sigma,\phi_\sigma)$ satisfies~(cf.~\eqref{regtheta}--\eqref{regBetaphi})\begin{align}
  & \theta_\sigma \in \H1H \cap \L\infty{\VA r} \cap \L2{\VA{2r}},
  \qquad
  \label{regthetasig}
  \\
  & \phi_\sigma \in \H1H \cap \L\infty{\VB\sigma} ,
  \label{regphisig}
  \\
  & \Beta(\phi_\sigma) \in \LQ1 .
  \label{regBetaphisig}
\end{align}
The convergence theorem we prove is as follows.
\Bthm
\label{Pier}
Under the assumptions \accorpa{hpconst}{defbetapi} and \accorpa{a-pier0}{a-pier4}, 
the family of solutions $(\theta_\sigma,\phi_\sigma)$ to the problem \accorpa{primasig}{regBetaphisig}
satisfies 
\Bsist
  & {\theta_{\sigma}} \to \theta
  & \quad \hbox{weakly star in $\H1H\cap\L\infty{\VA r}\cap \L2{\VA {2r}}$}\qquad 
  \non
  \\
  && \quad \hbox{and strongly in $\C0H$}\pier{,}
  \label{convthetasig}
  \\[2mm]
  & \phi_\sigma \to \phi
  & \quad \hbox{weakly in $\H1H$} 
  \label{convphisig}
\Esist
as $\sigma\searrow 0$, where the limit pair $(\theta,\phi)$ is the unique solution to the problem 
\Bsist
  && \dt\theta + \ell \, \dt\phi + A^{2r} \theta = f
  \quad \aeQ \, ,
  \label{prima-fin}
  \\[1mm]
  && \dt\phi + \phi - P\phi + \xi - \gamma\, \phi = \ell \, \theta , 
  \quad \hbox{for some \gianni{$\xi\in\L2H$ satisfying}}
  \non 
  \\
  && \quad
  \gianni{\xi \in \beta (\phi) \quad \aeQ \, ,  }
  \label{seconda-fin}
  \\[1mm]
  && \theta(0) = \thetaz \, 
  , \quad
  \phi(0) = \phiz 
  \label{cauchy-fin}
\Esist
and, in \eqref{seconda-fin}, $P$ denotes the $H$-projection operator on the kernel of the operator $B$.
\Ethm 

\Brem
\label{RemOperP}
With reference to Remark~\ref{Remhpsimple}, let us point out that 
in the case whether $B$ is the Laplace 
operator $-\Delta$ with Neumann boundary conditions, the operator 
$P$ maps any element $v\in H$ into a constant function, 
which is proportional to the mean value of $v$ 
as the first eigenfunction \gianni{of $B$} is $\eta_1 = |\Omega|^{- 1/2} $. 
\Erem

\Brem
\label{Remspecialpi}
It turns out that Theorem~\ref{Pier} works under the special assumption 
\eqref{a-pier1} of a linear function $\pi$ that is the derivative of a concave quadratic function~$\Pi$. 
This is of course a special situation, but let us point out that all 
the three significant examples of potentials \eqref{regpot}--\eqref{obspot} 
contemplate exactly a {\it quadratic concave perturbation}, along with a convex and possibly singular or nonsmooth function. \Erem

The \gianni{whole} section is devoted to the proof of this theorem. About the uniqueness of the limit
$(\theta,\phi)$ and the related continuous dependence property with respect to the data 
$(f, \theta_0, \phi_0)$, it is not difficult to somehow reproduce the proof given in Section~\ref{UNIQUENESS} for the case of a constant $\ell $ and derive
an estimate similar to \eqref{percontdep}, which leads~to  

\Bprop
\label{CDpier}
Assume that \accorpa{hpconst}{defbetapi} and \accorpa{a-pier0}{a-pier1} are satisfied.
\gianni{Moreover, let the data $f$, $\theta_0$ and $\phi_0$ satisfy \accorpa{a-pier2}{a-pier4}
for some family $\{ f_{\sigma},\, \theta_{0,\sigma}, \, \phi_{0,\sigma} \}$.}
Then, there exists a unique solution
$$ (\theta,\phi) \in \bigl(\H1H\cap\L\infty{\VA r} \cap \L2{\VA {2r}}\bigr) \times \H1H $$
to the problem \accorpa{prima-fin}{cauchy-fin}. 
Moreover, if $(f_i,\theta_{0i},\phi_{0i})$, $i=1,2$, are two choices of the data 
and $(\theta_i,\phi_i)$ are the corresponding solutions, then
we have
\begin{align}
  & \norma{\theta_1-\theta_2}_{\L2H}
  + \norma{1*(\theta_1-\theta_2)}_{\L\infty{\VA r}}
  \non\\
  &+ \norma{\phi_1-\phi_2}_{\L\infty H}
  + \norma{\phi_1 - P\phi_1 - \phi_2+ P \phi_2 }_{\L2{H}}  
  \non
  \\
  & \leq K \bigl(
    \norma{1*(f_1-f_2)}_{\L2H}
    + \norma{\theta_{01}-\theta_{02}}
    + \norma{\phi_{01}-\phi_{02}}
  \bigr)
  \label{contdeppier}
\end{align}
for some constant $K$ depending only on $\ell, \, \gamma , \, T$.
\Eprop

Now, we concentrate our efforts on the proof of the remaining convergence properties 
stated in Theorem~\ref{Pier}. We start by proving the following auxiliary result. 
\Blem
\label{Conv-lemma}
Assume that $v\in \VB{\sigma_0} $ for some $\sigma_0 >0$. Then we have that
$ B^\sigma v $ is well defined for all $\sigma \in [0,\sigma_0]$ and 
\Beq 
\label{pier1} 
B^\sigma v \to v - Pv \quad \hbox{strongly in } \, H \quad \hbox{as }\, \sigma\searrow 0,
\Eeq 
where, as above, $P$ denotes the $H$-projection on $\{ v\in D(B) : \ Bv = 0\}$. 
\Elem

\Bdim
The first part of the statement follows easily from \eqref{defdomBs} and \eqref{defArBs}. 
In particular,  we note that 
$$
B^\sigma v = \somma j1\infty \mu_j^\sigma (v,\eta_j) \eta_j \quad \hbox{for all } 0<\sigma \leq \sigma_0, \quad \ v-Pv = \sum_{\mu_j >0} (v,\eta_j) \eta_j . 
$$
Then, we have to prove that 
$$
B^\sigma v  - (v-Pv) = \sum_{\mu_j >0}  \left(\mu_j^\sigma -1 \right) (v,\eta_j) \eta_j 
\to 0\quad \hbox{strongly in } \, H \quad \hbox{as }\, \sigma\searrow 0.
$$
In view of \eqref{eigen}--\eqref{complete}, it is sufficient to verify that 
$$
\norma{B^\sigma v  - (v-Pv)}^2 = \sum_{\mu_j >0}  \left(\mu_j^\sigma -1 \right)^2 |(v,\eta_j)|^2 
$$
tends to $0$ as $\sigma \searrow 0$. We observe that $ \left(\mu_j^\sigma -1 \right)^2 \leq 
1$ if $\mu_j \leq 1$ and $ \left(\mu_j^\sigma -1 \right)^2 \leq 
 \mu_j^{2\sigma_0}$ if $\mu_j > 1$. Hence, we have that 
$$
\sum_{\mu_j >0}  \left(\mu_j^\sigma -1 \right)^2 |(v,\eta_j)|^2 \leq 
\somma j1\infty  \big(1 +\mu_j^{2\sigma_0} \big) |(v,\eta_j)|^2
= \norma{v}^2 + \norma{B^{\sigma_0} v}^2 < +\infty. 
$$ 
Therefore, the reader can realize that it is possible to apply the Lebesgue dominated 
convergence theorem, with respect to the counting measure~$\#$, to the family of functions 
$$ f_\sigma (j) = 
\begin{cases}
0 &  \hbox{if } \mu_j=0  \\[1mm]
 (\mu_j^\sigma -1 )^2 \gianni{\,|(v,\eta_j)|^2} & \hbox{if } \mu_j>0 
 \end{cases}
 \ , \quad  j\in \enne. 	
$$
\gianni{Since
\begin{align*}
&f_\sigma  \to 0 \quad \hbox{pointwise in }\, \enne \quad \hbox{as } \, \sigma\searrow 0 \,,
\\
&0 \leq f_\sigma (j) \leq g (j):= \big(1 +\mu_j^{2\sigma_0} \big) \, |(v,\eta_j)|^2 \quad \hbox{for all }\, j\in \enne \,,
\end{align*}
and $g$ is summable with respect to $\#$ by the above calculation,}
we can conclude that 
$$
\norma{B^\sigma v  - (v-Pv)}^2 = \int_\enne f_\sigma (j)\,d\# \,(j) \to 0
$$
and \eqref{pier1} and the lemma are completely proved.
\Edim

Next, we recall and take advantage of the uniform estimates pointed out in Section~\ref{EXISTENCE}. Arguing as in the derivation of \eqref{primastima}, 
recalling \eqref{defnormaBs} and using the 
lower semicontinuity properties when passing 
to the limit as $\eps \searrow 0$, from \eqref{perprimastima} we obtain
\begin{align}
  & \frac 12 \, \norma{\theta_\sigma(t)}^2
  + \iot \norma{A^r\theta_\sigma(s)}^2 \, ds
  + \frac 12 \intQt |\dt\phi_\sigma|^2
  \non\\
  &+ \frac 12 \left( \norma{\phi_\sigma(t)}^2 + 
 \norma{B^\sigma \phi_\sigma(t)}^2 \right) 
   + \iO \Beta(\phi_\sigma(t))
  \non \\ 
  &\leq c+ \frac 12 \,  \norma{\theta_{0,\sigma}}^2
  + \frac 12 \left( \norma{\phi_{0,\sigma}}^2 +
  \norma{B^\sigma\phi_{0,\sigma}}^2 \right)\non \\
  &\quad   + \iO \Beta({\phi_{0,\sigma}}) 
  + \frac 12\int_0^t \gianni{\norma{f_\sigma(s) + (1+\gamma)\phi_\sigma(s)}^2 \, ds} \,.
  \non
\end{align}
Hence, by virtue of \eqref{a-pier2}--\eqref{a-pier4} 
and applying the Gronwall lemma, we deduce that
\begin{align}
  &\norma{\theta_\sigma}_{\L\infty H\cap\L2{\VA r}}
  + \norma{\phi_\sigma}_{\H1H}
  \non\\
  &+ \norma{B^\sigma \phi_\sigma}_{\L\infty{H}}
      + \norma{\Beta (\phi_\sigma)}_{\L\infty\Luno}
  \leq c \,.
  \label{primastimasig}
\end{align}
Now, we can test \eqref{primasig} by $\dt\theta_\sigma$ and integrate with
respect to time obtaining
\Beq
  \intQt |\dt\theta_\sigma|^2
  + \frac 12 \, \norma{A^r\theta_\sigma(t)}^2
  =  \frac 12 \, \norma{A^r\theta_{0,\sigma}}^2  +
  \intQt \bigl( f_\sigma -\ell\,\dt\phi_\sigma \bigr) \dt\theta_\sigma \,.
  \non 
\Eeq
Then, the Young inequality, \eqref{primastimasig} and \eqref{a-pier2}--\eqref{a-pier3}
enable us to infer that
\Beq
  \norma{\theta_\sigma}_{\H1H\cap\L\infty{\VA r}}
  \leq c \,.
  \label{secondastimasig}
\Eeq
Moreover, in view of \eqref{a-pier2} and by a comparison in \gianni{\eqref{primasig}} we have that 
\Beq
  \norma{\theta_\sigma}_{\L2{\VA {2r}}}
  \leq c \,.
  \label{terzastimasig}
\Eeq
Due to \accorpa{primastimasig}{terzastimasig}, we can pass to 
the limit as $\sigma \searrow 0$, at the 
beginning for a subsequence, by using standard weak  
and known strong compactness results (see, e.g., \cite[Sect.~8, Cor.~4]{Simon}). Hence,
we find out the \gianni{convergence} \eqref{convthetasig}--\eqref{convphisig} to $\theta$ and $\phi$,
along with
\Beq  B^\sigma \phi_\sigma \to \zeta 
  \quad \hbox{weakly star in $\L\infty H$}.
  \label{convzetasig}
\Eeq
In a first verification, inspired by Lemma~\ref{Conv-lemma}, we aim to check the weak star limit $\zeta $ in \eqref{convzetasig} is nothing but $\phi - P\phi$. Actually, we show that 
\Beq B^\sigma \phi_\sigma \to \phi - P\phi
\quad \hbox{weakly in $\L2H$},
  \label{pier5}
\Eeq 
by verifying this property with respect to a dense subset of $\L2H$. Indeed, thanks to \eqref{complete} 
\piernew{it suffices to} prove that 
\begin{align}
&\int_0^T ( B^\sigma \phi_\sigma (t) , \psi (t) \eta_j) \, dt \to 
\int_0^T ( \phi (t) - P\phi (t) , \psi (t) \eta_j) \, dt \non \\
&\quad \hbox{ for all } \, \psi \in L^2(0,T) \, \hbox{ and } \, j\in \enne . 
\label{pier2} 
\end{align} 
Note that the integrals in \eqref{pier2} all vanish if the index $j $ is such that 
the eigenvalue $\mu_j $ is equal to $0$. If instead $\mu_j >0$, then with the help of \eqref{defArBs} and \eqref{convphisig}
we have that  
$$
\int_0^T ( B^\sigma \phi_\sigma (t) , \psi (t) \eta_j) \, dt = 
\mu_j^\sigma \int_0^T ( \phi_\sigma (t) , \psi (t) \eta_j) \, dt
\to \int_0^T ( \phi (t) , \psi (t) \eta_j) \, dt
$$
as $\sigma \searrow 0$.  
\gianni{Moreover, $\eta_j$~is orthogonal to the kernel of~$B$.
All this} means that in both cases there is convergence to 
$\int_0^T ( \phi (t) - P\phi (t) , \psi (t) \eta_j) \, dt$ and \eqref{pier2} is ensured. 

At this point, it remains to prove that the limiting pair $(\theta,\phi)$ solves the system \eqref{prima-fin}--\eqref{cauchy-fin}. The initial conditions \eqref{cauchy-fin} hold true
by virtue of \eqref{convthetasig}--\eqref{convphisig}, \eqref{cauchysig} and 
\eqref{a-pier3}--\eqref{a-pier4}. In particular, note that \eqref{convphisig}
guarantees at least that 
\Beq  \phi_\sigma (t) \to \phi (t) 
  \quad \hbox{weakly in $H$, for all } \, t \in [0,T].
  \label{pier6}
\Eeq
On the other hand, recalling \eqref{a-pier2} and passing to the limit in \eqref{primasig}
we arrive at \eqref{prima-fin}. Next, we let the test function $ v$ in \eqref{secondasig} to depend also
on time, taking $v\in \L2 {\VB{\sigma_0} }$, and multiply the inequality by $e^{-2\gamma t}$, then integrating over $(0,T)$. We obtain
\begin{align}
& \int_0^T \bigl( e^{-\gamma t}(\dt\phi_\sigma(t)-  \gamma \, \phi_\sigma(t))  ,e^{-\gamma t}  (\phi_\sigma - v)(t) \bigr)\, dt 
\non\\
 &+ \int_0^T e^{-2\gamma t}\bigl( B^\sigma \phi_\sigma(t) , B^\sigma (\phi_\sigma- v)(t) \bigr) \, dt
%  \non
%  \\
%  & \quad {}
  +\int_0^T \!\!\iO e^{-2\gamma t} \Beta(\phi_\sigma(t))\, dt
  \non
  \\
  & \leq \int_0^T e^{-2\gamma t}\bigl( \ell\, \theta_\sigma(t) , (\phi_\sigma - v) (t)  \bigr)\, dt
  +\int_0^T \!\!\iO e^{-2\gamma t} \Beta(v(t))\, dt 
  \non
  \\
  & 
  \quad 
  \hbox{for every $v\in\L2{\VB{\sigma_0}}$} \, .
  \label{nuovasig}
\end{align}
Let us introduce the family of functions
$$ 
\rho_\sigma (t): =  e^{-\gamma t} \phi_\sigma(t), \quad t\in [0,T],
$$ 
and point out that (cf.~\eqref{convphisig} and~\eqref{pier6})
\begin{align}
  &\rho_\sigma \to \rho \quad \hbox{weakly in $\H1H$}, \quad \hbox{where } \, \rho (t): =  e^{-\gamma t} \phi(t), \ \,  t\in [0,T],
  \label{pier3}\\
 &\rho_\sigma (t) \to \rho (t) 
  \quad \hbox{weakly in $H$, for all } \, t \in [0,T].
  \label{pier4}
\end{align} 
Then, we observe that 
\begin{align}
& \int_0^T \bigl( e^{-\gamma t}(\dt\phi_\sigma(t)-  \gamma \, \phi_\sigma(t))  ,e^{-\gamma t}  (\phi_\sigma - v)(t) \bigr)\, dt 
\non\\
 &
 =  \int_0^T \bigl( \dt\rho_\sigma(t)  ,\rho_\sigma (t) - e^{-\gamma t} v(t) \bigr) \, dt 
\non\\
 &= \frac12 \bigl\Vert \rho_\sigma(T)  \bigr\Vert^2  -  \frac12 \bigl\Vert \phi _{0,\sigma}  \bigr\Vert^2
- \int_0^T \bigl( \dt \rho_\sigma(t)  ,e^{-\gamma t} v(t) \bigr) \, dt \non
\end{align}
and, on account of \eqref{a-pier4}, \eqref{pier3}, \eqref{pier4} and the weak lower semicontinuity of norms, we have that 
\gianni{%
\begin{align}
  & \int_0^T \bigl( e^{-\gamma t}(\dt\phi(t)-  \gamma \, \phi(t))  ,e^{-\gamma t}  (\phi - v)(t) \bigr)\, dt 
  \non
  \\
  & = \int_0^T \bigl( \dt\rho(t)  ,\rho (t) - e^{-\gamma t} v(t) \bigr) \, dt  
  \non
  \\
  & = \frac12 \bigl\Vert \rho (T)  \bigr\Vert^2  -  \frac12 \bigl\Vert \phi _{0}  \bigr\Vert^2
  - \int_0^T \bigl( \dt \rho (t)  ,e^{-\gamma t} v(t) \bigr) \, dt  
  \non
  \\
  & \leq \liminf_{\sigma\searrow 0} \Bigl\{
    \frac12 \bigl\Vert \rho_\sigma(T)  \bigr\Vert^2  -  \frac12 \bigl\Vert \phi _{0,\sigma}  \bigr\Vert^2
    - \int_0^T \bigl( \dt \rho_\sigma(t)  ,e^{-\gamma t} v(t) \bigr) \, dt
  \Bigr\}
  \non
  \\
  & = \liminf_{\sigma\searrow 0} \int_0^T \bigl( e^{-\gamma t}(\dt\phi_\sigma(t) - \gamma \, \phi_\sigma(t))  ,e^{-\gamma t}  (\phi_\sigma - v)(t) \bigr)\, dt \,.
  \label{pier7}
\end{align}
}%
Similarly, we recall \eqref{pier5} and point out that  $ B^\sigma v \to v-Pv$ strongly in $\L2H $: 
indeed, this is a consequence of Lemma~\ref{Conv-lemma}, 
\gianni{the bounds
\Bsist
  && \norma{B^\sigma v(t)}
  \leq \norma{v(t)}_{B,\sigma}
  \leq c \, \norma{v(t)}_{B,\sigma_0}
  \quad \hbox{\aat\ and every $\sigma\in(0,\sigma_0]$\piernew{,}}
  \non
  \\
  &&\quad  \piernew{\hbox{along with }} \ \ioT \norma{v(t)}_{B,\sigma_0}^2 \, dt = \norma v_{\L2{\VB{\sigma_0}}}^2 < +\infty \, \piernew{,}
  \non
\Esist
}%
and the Lebesgue dominated convergence theorem. Hence, as $\phi - P \phi$ is orthogonal to every element of the kernel of $B$, we infer that 
\gianni{%
\begin{align}
& \int_0^T e^{-2\gamma t}\bigl( (\phi - P \phi) (t) ,  (\phi   - v ) (t) \bigr) \, dt
 \non \\
& = \int_0^T e^{-2\gamma t}\bigl( (\phi - P \phi) (t) ,  (\phi - P \phi) (t)  - (v - P v) (t) \bigr) \, dt
 \non \\
& \leq \liminf_{\sigma\searrow 0}
\int_0^T e^{-2\gamma t}\bigl( B^\sigma \phi_\sigma(t) , B^\sigma (\phi_\sigma- v)(t) \bigr) \, dt \,.
 \label{pier8}
\end{align}
}%
Next, we observe that the function
$$
v\mapsto \int_0^T \!\!\iO e^{-2\gamma t} \Beta(v(t))\, dt
$$
is convex and lower semicontinuous in $\L2H$, as one can easily verify. Therefore, since  $\phi_\sigma$ weakly converges to 
$\phi$ in  $\L2H$ (see \eqref{convphisig}), we have that
\Beq
\int_0^T \!\!\iO e^{-2\gamma t} \Beta(\phi )\, dt  \leq\liminf_{\sigma\searrow 0} \int_0^T \!\!\iO e^{-2\gamma t} \Beta(\phi_\sigma (t))\, dt . 
\label{pier9} 
\Eeq
Now, we take advantage of \eqref{pier7}--\eqref{pier9} and, in view of  
\eqref{convthetasig} and \eqref{convphisig}, from \eqref{nuovasig} we deduce that 
\begin{align}
&\int_0^T \bigl( e^{-\gamma t}(\dt\phi(t)-  \gamma \, \phi(t))  ,e^{-\gamma t}  (\phi - v)(t) \bigr)\, dt
\non \\
& + \int_0^T e^{-2\gamma t}\bigl( (\phi - P \phi) (t) ,  (\phi   - v ) (t) \bigr) \, dt
+ \int_0^T \!\!\iO e^{-2\gamma t} \Beta(\phi )\, dt  
\non \\
& \leq\liminf_{\sigma\searrow 0}
\left(\int_0^T \bigl( e^{-\gamma t}(\dt\phi_\sigma(t)-  \gamma \, \phi_\sigma(t))  ,e^{-\gamma t}  (\phi_\sigma - v)(t) \bigr)\, dt \right.
\non\\
&\qquad\qquad\quad  \left.+ \int_0^T e^{-2\gamma t}\bigl( B^\sigma \phi_\sigma(t) , B^\sigma (\phi_\sigma- v)(t) \bigr) \, dt 
  +\int_0^T \!\!\iO e^{-2\gamma t} \Beta(\phi_\sigma(t))\, dt  \right)  \non
  \\
& \leq \int_0^T e^{-2\gamma t}\bigl( \ell\, \theta(t) , (\phi - v) (t)  \bigr)\, dt
  +\int_0^T \!\!\iO e^{-2\gamma t} \Beta(v(t))\, dt 
  \non
  \\
  & 
  \quad 
  \hbox{for every $v\in\L2{\VB{\sigma_0}}$} \, . \non
\end{align}
Therefore, at the end we derive the same inequality as \eqref{nuovasig} 
for the limit functions $\theta$ \gianni{and~$\phi$}. 
Moreover, it is not difficult to check that 
this inequality is equivalent to 
\begin{align} 
  & \bigl( \dt\phi(t) , \phi(t) - v \bigr)
  + \bigl( \phi(t) - P \phi (t) ,  \phi (t)   - v \bigr) 
%  \non
%  \\
%  & \quad {}
  + \iO \Beta(\phi(t))
  -\bigl( \gamma \, \phi(t) , \phi(t) - v \bigr)
  \non
  \\
  & \leq \bigl( \ell\, \theta(t) , \phi(t) - v \bigr)
  + \iO \Beta(v)
%  \non
%  \\
%  & 
  \quad 
  \hbox{\aat\ and every $v\in\VB{\sigma_0}$} \, .
  \label{pier10}
\end{align}
Then, by a density argument it is straightforward to infer that \eqref{pier10} holds true 
for all $v\in H$, whence the definition of subdifferential for the convex functional 
$$ v \in H \mapsto  \iO \Beta(v) \in [0,+\infty] $$ 
enables us to conclude that \gianni{\aat
\Beq
\xi(t) := (\ell \, \theta - \dt\phi - \phi + P  \phi + \gamma\, \phi)(t)  \in \beta (\phi(t) )   
\quad \aeO \,.
\label{pier11}
\Eeq
}%
It is easy now to see that \eqref{pier11} finally leads to \eqref{seconda-fin}.

%%%%%%%%%%%%%%%%%%%%%%%%%%%%%%%%%%%%%%%%%%%%%%%%%%%%%%%%%%%%%%%%%%%%%%%%

\section*{Acknowledgments}
\pier{This research was supported by the Italian Ministry of Education, 
University and Research (MIUR): Dipartimenti di Eccellenza Program (2018�2022) 
-- Dept. of Mathematics �F. Casorati", University of Pavia. 
In addition, the authors gratefully acknowledge some other 
financial support from the MIUR-PRIN Grant 2015PA5MP7 ``Calculus of Variations'',}
the GNAMPA (Gruppo Nazionale per l'Analisi Matematica, 
la Probabilit\`a e le loro Applicazioni) of INdAM (Isti\-tuto 
Nazionale di Alta Matematica) and the IMATI -- C.N.R. Pavia. 

%%%%%%%%%%%%%%%%%%%%%%%%%%%%%%%%%
%% bibliography
%%%%%%%%%%%%%%%%%%%%%%%%%%%%%%%%%

\vspace{3truemm}

\Begin{thebibliography}{10}

{%
\bibitem{AD}
N. Abatangelo, L. Dupaigne,
Nonhomogeneous boundary conditions for the spectral fractional Laplacian,
Ann. Inst. H. Poincar\'e Anal. Non Lin\'eaire {\bf 34} (2017), 439-467.
\bibitem{AM}
M. Ainsworth, Z. Mao, Analysis and approximation of a fractional
Cahn--Hilliard equation, SIAM J. Numer. Anal. {\bf 55} (2017), 1689-1718.
\bibitem{AkSS1}
G. Akagi, G. Schimperna, A. Segatti,
Fractional Cahn--Hilliard, Allen--Cahn and porous medium equations,
J. Differential Equations {\bf 261} (2016), 2935-2985.
\bibitem{AkSS2}
G. Akagi, G. Schimperna, A. Segatti,
Convergence of solutions for the fractional Cahn--Hilliard system,
preprint arXiv:1801.01722 [math.AP] (2018), pp.~1-43.}
\bibitem{Barbu}
V. Barbu,
``Nonlinear Differential Equations of Monotone Type in Banach Spaces'',
Springer,
London, New York, 2010.

\bibitem{BFV}
M. Bonforte, A. Figalli, J.L. V\'azquez, 
Sharp global estimates for local and nonlocal porous medium-type equations in bounded 
domains, Anal. PDE {\bf 11} (2018), 945-982. 
\bibitem{BSV}
M. Bonforte, Y. Sire, J.L. V\'azquez,  
Existence, uniqueness and asymptotic behaviour for fractional porous medium on bounded domains, 
Discrete Contin. Dyn. Syst. {\bf 35} (2015), 5725-5767.
\bibitem{BV}
M. Bonforte, J.L. V\'azquez,  
A priori estimates for fractional nonlinear 
degenerate diffusion equations on bounded domains, 
Arch. Ration. Mech. Anal. {\bf 218} (2015), 317-362.

\bibitem{BSY}
L. Brasco, M. Squassina, Y. Yang,
Global compactness results for nonlocal problems,
Discrete Contin. Dyn. Syst. Ser. S {\bf 11} (2018), 391-424.

\bibitem{Brezis}
H. Brezis,
``Op\'erateurs maximaux monotones et semi-groupes de contractions
dans les espaces de Hilbert'',
North-Holland Math. Stud.
{\bf 5},
North-Holland,
Amsterdam,
1973.

\bibitem{CT}
X. {Cabr\'e},  J. Tan, 
Positive solutions of nonlinear problems involving the square root of the 
{L}aplacian, Adv. Math. {\bf 224} (2010), 2052-2093.
\bibitem{CS}
L.A. Caffarelli, P.R. Stinga,
Fractional elliptic equations, {C}accioppoli estimates and regularity,
Ann. Inst. H. Poincar\'e Anal. Non Lin\'eaire {\bf 33} (2016), 767-807.
\bibitem{caginalp}
G. Caginalp,
{An analysis of a phase field model of a free boundary},
Arch. Rational Mech. Anal.
{\bf 92}
(1986)
205-245.

\bibitem{cagnish}
G. Caginalp and Y. Nishiura,           
The existence of travelling waves for phase field 
equations and convergence to sharp interface models in the singular limit,
Quart. Appl. Math. {\bf 49} (1991) 147-162.

\bibitem{CCFT}
M. Caputo, M.  Ciarletta, M. Fabrizio, V. Tibullo,  
Melting and solidification of pure metals by a phase-field model.
Atti Accad. Naz. Lincei Rend. Lincei Mat. Appl.  {\bf 28}  (2017), 463-478.
\bibitem{CGS6}
P. Colli, G. Gilardi, J. Sprekels,
On an application of Tikhonov's fixed point theorem
to a nonlocal Cahn--Hilliard type system modeling phase separation,
J. Differential Equations {\bf 260} (2016), 7940-7964.
\bibitem{CGS18}
P. Colli, G. Gilardi, J. Sprekels,
Well-posedness and regularity for a generalized fractional Cahn--Hilliard system,
{preprint arXiv:1804.11290 [math.AP] (2018), pp.~1-35.}

\bibitem{g}
   H.~Gajewski,
   On a nonlocal model of non-isothermal phase separation,
   Adv. Math. Sci. Appl.,
   {\bf 12} (2002),
   569-586.

\bibitem{gz}
   H.~Gajewski and K.~Zacharias,
   On a nolocal~phase~separation~model,
   J. Math. Anal. Appl., {\bf 286} (2003),
   11-31.
\bibitem{GalEJAM}
C.G. Gal, 
Non-local Cahn--Hilliard equations with fractional dynamic boundary conditions, European J. Appl. Math. {\bf 28} (2017), 736-788. 
\bibitem{GalDCDS}
C.G. Gal, 
On the strong-to-strong interaction case 
for doubly nonlocal Cahn--Hilliard equations,
Discrete Contin. Dyn. Syst.  {\bf 37} (2017),  131-167.
\bibitem{GalAIHP}
C.G. Gal, 
Doubly nonlocal Cahn--Hilliard equations, 
Ann. Inst. H. Poincar\'e Anal. Non Lin\'eaire {\bf 35} (2018), 357-392.

\bibitem{Gru1} 
G. Grubb,  
Regularity of spectral fractional Dirichlet and Neumann problems,
Math. Nachr.  {\bf 289} (2016), 831-844.
\bibitem{Gru2} 
G. Grubb,  
Regularity in $L_p$ Sobolev spaces of solutions to fractional heat equations,
J. Funct. Anal.  {\bf 274}  (2018), 2634-2660.

\bibitem{HTY}
T. Hou, T. Tang, J. Yang, Numerical analysis of fully discretized Crank-Nicolson scheme for fractional-in-space Allen-Cahn equations. J. Sci. Comput. {\bf 72}  (2017), 1214-1231.

\bibitem{Kwa}
M. Kwa\'snicki, Ten equivalent definitions of the fractional Laplace operator,
Fract. Calc. Appl. Anal. {\bf 20} (2017), 7-51.

\bibitem{LWY}
Z. Li, H. Wang, D.  Yang, A space-time fractional phase-field model with tunable sharpness and decay behavior and its efficient numerical simulation, J. Comput. Phys.  {\bf 347}  (2017), 20-38.

{%
\bibitem{MN}
R. Musina, A.I. Nazarov, 
Variational inequalities for the spectral fractional Laplacian,
Comput. Math. Math. Phys.  {\bf 57} (2017), 373-386.
\bibitem{RS1}
X. Ros-Oton, J. Serra, 
The Dirichlet problem for the fractional Laplacian: regularity up to
the boundary, J. Math. Pures Appl. (9) {\bf 101} (2014), 275-302.
\bibitem{RS2}
X. Ros-Oton, J. Serra, 
The extremal solution for the fractional Laplacian, 
Calc. Var. Partial Differential Equations {\bf 50} (2014), 723-750.
\bibitem{SV0}
R. Servadei, E. Valdinoci, 
Variational methods for non-local operators of elliptic type, Discrete
Contin. Dyn. Syst. {\bf 33} (2013), 2105-2137.
\bibitem{SV1}		
R. Servadei, E. Valdinoci,  
On the spectrum of two different fractional operators,
Proc. Roy. Soc. Edinburgh Sect. A  {\bf 144}  (2014), 831-855.
\bibitem{SV2}		
R. Servadei, E. Valdinoci,  
Weak and viscosity solutions of the fractional Laplace equation,
Publ. Mat.  {\bf 58}  (2014), 133-154.
\bibitem{SV3}		
R. Servadei, E. Valdinoci, 
The Brezis-Nirenberg result for the fractional Laplacian,
Trans. Amer. Math. Soc.  {\bf 367} (2015), 67-102.}

\bibitem{Simon}
J. Simon,
Compact sets in the space $L^p(0,T; B)$,
Ann. Mat. Pura Appl.~(4) {\bf 146}, (1987), 65-96.

\End{thebibliography}

\End{document}

%%%%%%%%%%%%%%%%%%%%%%%%%%%%%%%%%%%%%%%%%%%%%%%%